\documentclass[12pt]{amsart}
\usepackage[dvips]{graphicx}
\newtheorem{thm}{Theorem}[section]
\newtheorem{cor}[thm]{Corollary}
\newtheorem{lem}[thm]{Lemma}
\newtheorem{prop}[thm]{Proposition}
\newtheorem{defn}[thm]{Definition}
\newtheorem{note}[thm]{Note}
\newtheorem{ex}[thm]{Example}

\newtheorem{corollary}[thm]{Corollary}

\begin{document}

\title[Finite Reflection Groups]
{Lattices in Finite Real Reflection Groups}
\author[T.~Brady]{Thomas Brady}
\address{School of Mathematical Sciences\\
Dublin City University\\
Glasnevin, Dublin 9\\
Ireland}
\email{tom.brady@dcu.ie}
\author[C.~Watt]{Colum Watt}
\address{School of Mathematical Sciences\\
Dublin Institute of Technology\\
Kevin St., Dublin 8\\
Ireland}
\email{colum@maths1.kst.dit.ie}
%\subjclass{}
%\keywords{}
\maketitle
\textbf{Abstract:}
For a finite real reflection group $W$ with Coxeter element $\gamma$
we give a uniform proof that the closed interval, $[I, \gamma]$
forms a lattice in the partial order on $W$ induced by reflection length.
The proof involves the construction of a simplicial complex which can be
embedded in the type W simplicial generalised associahedron.

\vskip .2cm
\section{Introduction.}\label{s:intro}
Let $W$ be a finite real reflection group. Associated to $W$ is a
finite type Artin group or generalised braid group, $A(W)$. Much of the
work on finite type Artin groups takes Garside's paper, \cite{garside},
as its starting point, using the set of fundamental reflections as a
generating set for $W$ and a corresponding standard generating set
for $A(W)$. Recently it has been shown that when the set of {\em all}
reflections is used as a generating set for $W$ and a corresponding
generating set is used for $A(W)$, a parallel theory can be constructed.
In particular, the new positive monoid embeds in $A(W)$ and new
K$(\pi, 1)$'s for $A(W)$ have been constructed.  The larger generating
set gives $A(W)$ a second structure as a Garside group.  The larger
generating set is proposed in \cite{brady3}, but Daan Krammer
had used it independently in unpublished work. In the case of the braid
group $B_n$, where $W$ is the symmetric group $\Sigma_n$, the larger
generating set coincides with the band generators from \cite{bkl}.
The second Garside structure for $B_n$ is described in \cite{bdm}
and the structure for general finite $W$ in \cite{bessis}. The general
construction of $K(A(W),1)$'s is described in \cite{cmw}, using ideas from
\cite{bestvina}.
\vskip .2cm
A central result needed in the development of this parallel theory is
that the closed interval, $[I, \gamma]$, bounded by the identity $I$
and a Coxeter element $\gamma$, forms a lattice in the partial order
on $W$ induced by reflection length. The lattice property is used to
prove the embedding of the positive monoid and the asphericity of the
new K$(\pi, 1)$'s.  Existing proofs of the lattice property use the
classification of finite real reflection groups with different methods
applied to the different groups. The symmetric group is handled in
\cite{bradysym}, the $C_n$ and $D_n$ groups in
\cite{bwcd}. Bessis treats all cases in \cite{bessis}.
\vskip .2cm
In this paper we give a new proof of the lattice property that is
independent of the classification of finite real reflection groups.
For this purpose we introduce a new simplicial complex which is a
geometric model for the partially ordered set $[I,\gamma]$.
If $n$ is the rank of $W$, this simplicial complex will lie in the
$(n-1)$-sphere, $S^{n-1}$, in ${\bf R}^n$ and its vertex set will consist
of a set of positive roots for $W$.
\vskip .2cm
The layout of the paper is as follows. In section~\ref{s:poset}, we recall
definitions and results about the partial order and spherical
simplicial complexes. In section~\ref{s:root-mu} we extend some of the
material from \cite{steinberg} and prove some properties of the dot
products of roots of $W$ with vertices of a Petrie polygon. In
section~\ref{s:Xgamma} we define the simplicial complex $X(\gamma)$ and
derive some of its properties. In section~\ref{s:subcomplex}  we define the subcomplexes
$X(\sigma)$ for $\sigma \le \gamma$. After some preparatory work in section~\ref{s:walls},
we characterise the geometric realisation of $X(\sigma)$ and we prove the lattice property in
section~\ref{s:lattice}. We conclude in section~\ref{s:fomin} by explaining the
connection between our construction and the generalised associahedra of
\cite{fz}. In the process, we give a new proof that the latter is a spherical
simplicial complex.

\section{The partial order.} \label{s:poset}
Let $W$ be a finite, irreducible real reflection group with
reflection set ${\mathcal R}$. We recall from \cite{bwcd} the
partial order on $W$ given by the reflection length function~$l$.
For $u, w \in W$ we say
\[u \le w \quad \quad \Leftrightarrow \quad \quad l(w) = l(u) +
l(u^{-1}w),\]
where $l(v)$ is the smallest positive integer $k$ such that $v$ can be
written as a product of
$k$ reflections from ${\mathcal R}$.
Thus $u \le w$ if and only if there is a shortest
factorisation of $u$ as a product of reflections which is a prefix of a
shortest
factorisation of $w$.
\vskip .2cm
We note that the partial order on $W$ is the restriction to $W$ of a
partial order on the
orthogonal group $O(n)$ which is introduced in \cite{bworth}.
Some implications of this relationship are investigated in \cite{bwcd}.
Some of the notation and
results from those papers will be used here. In particular, if
$A \in O(n)$, we associate to $A$ two subspaces of ${\bf R}^n$, namely
\[M(A) = \mbox{im}(A-I) \quad \mbox{ and } \quad F(A) =
\mbox{ker}(A-I),\]
which we  call the moved space of $A$ and the fixed space of $A$
respectively.
We recall that $M(A) = F(A)^{\perp}$ and note that $F(A) = F(A^{-1})$. The main
result of \cite{bworth} implies that if $V$ is a subspace of $M(A)$ then
there is a
unique $B \in O(n)$ satisfying $M(B) = V$ and $B \le A$. It was noted in
\cite{bwcd} that a group element $\alpha \le \gamma$ is characterised by its reflection
set $S_{\alpha}$, where
\[S_{\alpha} = \{R \in {\mathcal R} \mid R \le \alpha\}.\]
It was also noted in \cite{bwcd} that whenever $\alpha, \beta $ and $\delta$
are elements of the interval~$[I,\gamma]$ in $W$ with the properties that
$\delta \le \alpha, \beta$ and
$ S_{\delta} = S_{\alpha} \cap S_{\beta}$,
then $\delta $ must be the greatest lower bound of $\alpha$ and $\beta$ in
$W$; that is,
if $\tau \in W$ satisfies $\tau \le \alpha, \beta$ then $\tau \le
\delta$.
\vskip .2cm
We record for convenience here some results which either appear in \cite{bworth}
and \cite{bwcd} or are elementary consequences.
\begin{equation}\label{e:posetmap}
\mbox{If } \alpha \le \beta \mbox{ then } M(\alpha) \subset M(\beta) \mbox{ and }
F(\beta) \subset F(\alpha).
\end{equation}
If $R(\vec v)$ denotes the reflection in $\vec v^{\perp}$ then
\begin{equation}\label{e:conj}
\vec v \in M(\alpha) \Rightarrow R(\vec v)\alpha = \alpha R[\alpha^{-1}(\vec v)].
\end{equation}
\begin{equation}\label{e:knockoff}
\mbox{If } \alpha \le \beta \le \delta \mbox{ then }
\alpha^{-1}\beta \le \alpha^{-1}\delta \mbox{ and  }
\beta\alpha^{-1} \le \delta\alpha^{-1}.
\end{equation}
\begin{equation}\label{e:induce}
\mbox{If } \alpha , \beta \le \delta \mbox{ and  }
M(\alpha) \subset M(\beta) \mbox{ then } \alpha \le \beta.
\end{equation}
If $l(\alpha) = k$ and $R_i$ are reflections then
\begin{equation}\label{e:rhoperpmu}
\alpha = R_1\dots R_k \Rightarrow M(R_j) \subset M(R_1\alpha)
\mbox{ for } j = 2, \dots k.
\end{equation}
\begin{equation}\label{e:algpo}
\alpha \le \beta \le \delta \Rightarrow \beta^{-1}\delta \le \alpha^{-1}\delta
\mbox{ and  } \delta\beta^{-1} \le \delta\alpha^{-1}.
\end{equation}
If $R_1$ and $R_2$ are distinct reflections then
\begin{equation}\label{e:calc}
R_1R_2 \le \alpha \Leftrightarrow R_2 \le R_1\alpha
\Leftrightarrow R_1 \le \alpha R_2.
\end{equation}
\vskip .2cm
The complexes that we define in sections~\ref{s:Xgamma} and \ref{s:subcomplex}
are spherical simplicial complexes. We will
use (often without mention) the following facts about spherical simplices.
Any linearly independent set
of vectors in ${\bf R}^n$ determines a spherical simplex provided the
angle
between each pair of vectors is less than $\pi$. For such a set the
spherical simplex is obtained by
intersecting the unit sphere with the positive cone on those vectors.
The spherical simplex determined by the linearly
independent set
$\{v_1, v_2, \dots ,v_k\}$ will be denoted by $\langle v_1, v_2, \dots
,v_k \rangle$.
Since any fixed set of positive roots for $W$ lies in a single open
halfspace, each linearly independent subset of such a set of positive roots
determines a spherical simplex.

\section{Roots of $W$ and vertices on the Petrie polygon.}\label{s:root-mu}
In \cite{steinberg} Steinberg gives new proofs of several results about
irreducible, finite real reflection groups without using the classification
of these groups.  He uses the fact that, when the group is finite,
the fundamental reflections can be partitioned into two mutually
orthogonal subsets.  We will use his numbering of the roots of $W$ in
the construction of our simplicial complex.  Fix a fundamental chamber $C$ of the
standard tesselation
with inward unit normals $\alpha_1, \dots, \alpha_n$ and let
$R_1, \dots , R_n$ be the corresponding reflections. Since $W$
is finite we can assume that it is possible to choose the
ordering on the inward normals so that
$S_1 = \{\alpha_1, \dots , \alpha_s\}$ and $S_2 = \{\alpha_{s+1}, \dots , \alpha_n\}$
are orthonormal sets (from Lemma 2.2 of \cite{steinberg}). Let
$\gamma$ denote the Coxeter element given by
$\gamma = R_1 R_{2} \dots R_n$ and let $h$ denote the order of
$\gamma$ in $W$.  We note that $M(\gamma)$ is all of ${\bf R}^n$.  As in \cite{steinberg}, we set
$\rho_i = R_1R_2\dots R_{i-1}\alpha_i$,
where the $\alpha$'s and the $R$'s are indexed cyclically modulo $n$.
We will use the following explicit formulae which are
easily verified.
\[
\rho_i = \left\{
\begin{array}{cl}
\alpha_i &\mbox{ for } i = 1, \dots, s\\
-\gamma (\alpha_i)&\mbox{ for } i = s+1, \dots , n\\
\gamma(\rho_{i-n})&\mbox{ for } i > n.
\end{array}
\right.
\]
It is proved in \cite{steinberg} that the positive roots relative to
$C$ are $\rho_1, \dots , \rho_{nh/2}$ while  the negative roots are
$\rho_{(nh/2)+1}, \dots , \rho_{nh}$.  Furthermore the last
$n-s$ positive roots are a permutation (possibly trivial) of $S_2$.
For each $i$, let  $R(\rho_i)$
denote the reflection in $W$ with fixed subspace $\rho_i^{\perp}$
which is given by
\[R(\rho_i)(x) = x-2(\rho_i \cdot x)\rho_i.\]
\begin{note}\label{n:consec}
If $\rho_i, \rho_{i+1}, \dots , \rho_{i+n-1}$ are $n$ consecutive
roots in this ordering then
$R(\rho_i)R(\rho_{i+1})\dots R(\rho_{i+n-1}) = \gamma^{-1}$.
\end{note}
Let $\{\beta_1, \dots,
\beta_n\}$ be the dual basis to $\{\alpha_1, \dots, \alpha_n\}$.
Thus
\[\beta_i \cdot \alpha_j = \delta_{ij}, \mbox{ for } 1 \le i,j \le n,
\mbox {where } \delta_{ij} \mbox{ is the Kronecker delta.}\]
With a slight change of notation from \cite{steinberg} we define the vectors
$\mu_i$ for $1 \le i \le nh$ by
$\mu_i = R_1R_2\dots R_{i-1}\beta_i$,
where the $\beta$'s and the $R$'s are again indexed cyclically modulo $n$.
As with the $\rho$'s the following explicit formulae are readily verified.
\[\mu_i = \left\{
\begin{array}{ll}
\beta_i &\mbox{ for } i = 1, \dots, n\\
\gamma\mu_{i-n}&\mbox{ for } i > n.
\end{array}
\right.\]

\begin{thm}\label{t:petrieroot}
$\gamma \mu_i = \mu_i - 2\rho_i$ for $1 \le i \le nh$.
\end{thm}
\emph{Proof:}\hspace{.2cm}
Because of the recursions satisfied by both $\rho_i$ and $\mu_i$
it is sufficient to establish this result for $1 \le i
\le n$. In this case $\mu_i=\beta_i$ and, since $R_j(\beta_i)= \beta_i$ if
$j \ne i$, we have
\begin{eqnarray*}
\gamma \beta_i &=& (R_1\dots R_n)\beta_i\\
&=& (R_1\dots R_i)\beta_i\\
&=& (R_1\dots R_{i-1})(\beta_i-2\alpha_i)\\
&=& \beta_i-2(R_1\dots R_{i-1})\alpha_i\\
&=& \beta_i-2 \rho_i  \ \ \ \mbox{as required.}
\end{eqnarray*}
\begin{cor}\label{c:petrieroot}
With $\rho_i$ and $\mu_i$
defined as above we have
\begin{description}
\item[\rm (a)] $(\gamma - I) \mu_i = -2\rho_i$
\item[\rm (b)] $\mu_i \cdot \rho_i = 1$
\item[\rm (c)] $\mu_i \in F(R(\rho_i)\gamma)$
\end{description}
\end{cor}
\emph{Proof:}\hspace{.2cm}
Part (a) is immediate from theorem \ref{t:petrieroot}.
Part  (b) follows from theorem \ref{t:petrieroot}
and the fact that $\mu_i \cdot \mu_i = \gamma \mu_i \cdot \gamma
\mu_i$. From theorem~\ref{t:petrieroot} and part (b) we see that
$\gamma(\mu_i) = R(\rho_i)(\mu_i)$, which is equivalent to
statement (c).
\begin{note}
We observe that parts (b) and (c) of corollary~\ref{c:petrieroot}
characterise $\mu_i$ as the unique vector in the one-dimensional
subspace $F(R(\rho_i)\gamma)$ satisfying $\mu_i\cdot \rho_i = 1$.
\end{note}
\begin{note}
The vectors $\mu_i$ are shown in \cite{steinberg} to be
the vertices of a Petrie polygon. In particular the vectors
$\mu_1, \dots , \mu_n$ are
the vertices of the fundamental chamber $C$ while the other $\mu$'s are
obtained by
appropriate reflections using $\gamma(\mu_i) = R(\rho_i)(\mu_i)$.
It follows from theorem~\ref{t:petrieroot} that the  invertible linear transformation
$I- \gamma$ interchanges
the vertices of a Petrie polygon with the set of all roots; in so doing,
it takes the first $nh/2$
vertices to the set of positive roots.
\end{note}
In what follows we make extensive use the properties of the matrix of
dot products $[\mu_i\cdot \rho_j]$. Before establishing
the general result we will present an example.

\begin{ex}\label{e:icos}
Consider the symmetry group~$I_3$ of the icosohedron.
A simple system for $I_3$ is
\[\alpha_1 = (1,0,0), \alpha_2 = (0,1, 0),
\alpha_3 = (1/2)(-1, -\tau ,\tau -1)\]
where $\tau = 2\cos(\pi /5) = (1+\sqrt{5})/2$ is the positive root
of $x^2-x-1$. Thus $\tau$ satisfies
$\tau^2 = \tau +1$ and the other root of $x^2-x-1$ is
$-1/\tau =1-\tau$.

The dual basis to $\{\alpha_1,\alpha_2,\alpha_3\}$ is
\[\beta_1 = (1,0 ,\tau), \beta_2 = (0,1, \tau+1), \beta_3 =
(0, 0, 2\tau).\]
From the formula for a reflection we compute
\[R_1(\alpha_1) = -\alpha_1, \,\,\,R_1(\alpha_2) = \alpha_2,\,\,\,
R_1(\alpha_3) = \alpha_3+ \alpha_1.\]
\[R_2(\alpha_1) = \alpha_1, \,\,\,R_2(\alpha_2) = -\alpha_2,\,\,\,
R_2(\alpha_3) = \alpha_3 +\tau \alpha_2.\]
\[R_3(\alpha_1) = \alpha_1 + \alpha_3,\,\,\,
R_3(\alpha_2) = \alpha_2+ \tau \alpha_3,\,\,\,
R_3(\alpha_3) = -\alpha_3.\]
It will also be useful to know that
\begin{eqnarray*}
R_1R_2R_3(\alpha_1) &=& \tau \alpha_2 +\alpha_3\\
R_1R_2R_3(\alpha_2) &=& \tau \alpha_1 + \tau \alpha_2 + \tau \alpha_3\\
R_1R_2R_3(\alpha_3) &=& -\alpha_1 - \tau \alpha_2 -\alpha_3
\end{eqnarray*}
These values can be used to generate the
following table.
\vskip .2cm
\begin{tabular}{|l|l|l|l|}
\hline
$i$&$\rho_i$&$\mu_i$\\
\hline
$1$&$(1,0, 0)$&$(1,0,\tau)$\\
$2$&$(0, 1, 0)$&$(0,1, \tau+1)$\\
$3$&$(1/2)(1, \tau, \tau-1)$&$(0,0, 2\tau )$\\
$4$&$(1/2)(-1,\tau,\tau -1)$&$(-1,0, \tau)$\\
$5$&$(1/2)(\tau,\tau-1, 1)$&$(0, -1, \tau +1)$\\

$6$&$(1/2)(\tau -1, 1, \tau)$&$(-1,-\tau, \tau +1)$\\
$7$&$(1/2)(\tau, -\tau +1, 1)$&$(0,-\tau , 1 )$\\
$8$&$(1/2)(-\tau +1, 1, \tau )$&$(-\tau,-\tau, \tau)$\\
$9$&$(0, 0, 1)$&$(-\tau,-\tau -1, 1)$\\
$10$&$(1/2)(-\tau, \tau -1, 1)$&$(-\tau, -1, 0)$\\

$11$&$(1/2)(\tau -1,-1,\tau)$&$(-1,-\tau-1, 0)$\\
$12$&$(1/2)(-\tau +1,-1,\tau )$&$(-\tau,-\tau-1, -1)$\\
$13$&$(1/2)(1, -\tau ,\tau - 1)$&$(0,-\tau, -1)$\\
$14$&$(1/2)(-\tau , -\tau +1, 1)$&$(-\tau,-\tau, -\tau)$\\
$15$&$(1/2)(-1, -\tau , \tau -1)$&$(-1,-\tau, -\tau -1)$\\
\hline
\end{tabular}
\vskip .2cm
The entry in the $i^{\rm th}$ row and $j^{\rm th}$ column of the following table
is $\mu_i\cdot \rho_j$.
\vskip .2cm
\begin{tabular}{|r|r|r|r|r|r|r|r|r|r|r|r|r|r|r|}
\hline
$1$&$0$&$1$
&$0$&$\tau$&$\tau$
&$\tau$&$1$&$\tau$
&$0$&$\tau$&$1$
&$1$&$0$&$0$\\
$0$&$1$&$\tau$
&$\tau$&$\tau$&$\tau^2$
&$1$&$\tau^2$&$\tau^2$
&$\tau$&$\tau$&$\tau$
&$0$&$1$&$0$\\
$0$&$0$&$1$
&$1$&$\tau$&$\tau^2$
&$\tau$&$\tau^2$&$2\tau$
&$\tau$&$\tau ^2$&$\tau^2$
&$1$&$\tau$&$1$\\
$-1$&$0$&$0$
&$1$&$0$&$1$
&$0$&$\tau$&$\tau$
&$\tau$&$1$&$\tau$
&$0$&$\tau$&$1$\\
$0$&$-1$&$0$
&$0$&$1$&$\tau$
&$\tau$&$\tau$&$\tau^2$
&$1$&$\tau^2$&$\tau^2$
&$\tau$&$\tau$&$\tau$\\
$-1$&$-\tau$&$-1$
&$0$&$0$&$1$
&$1$&$\tau$&$\tau^2$
&$\tau$&$\tau^2$&$2\tau$
&$\tau$&$\tau ^2$&$\tau^2$\\
$0$&$-\tau$&$-1$
&$-1$&$0$&$0$
&$1$&$0$&$1$
&$0$&$\tau$&$\tau$
&$\tau$&$1$&$\tau$\\
$-\tau$&$-\tau$&$-\tau$
&$0$&$-1$&$0$
&$0$&$1$&$\tau$
&$\tau$&$\tau$&$\tau^2$
&$1$&$\tau^2$&$\tau^2$\\
$-\tau$&$-\tau^2$&$-\tau^2$
&$-1$&$-\tau$&$-1$
&$0$&$0$&$1$
&$1$&$\tau$&$\tau^2$
&$\tau$&$\tau^2$&$2\tau$\\
$-\tau$&$-1$&$-\tau$
&$0$&$-\tau$&$-1$
&$-1$&$0$&$0$
&$1$&$0$&$1$
&$0$&$\tau$&$\tau$\\
$-1$&$-\tau^2$&$-\tau^2$
&$-\tau$&$-\tau$&$-\tau$
&$0$&$-1$&$0$
&$0$&$1$&$\tau$
&$\tau$&$\tau$&$\tau^2$\\
$-\tau$&$-\tau^2$&$-2\tau$
&$-\tau$&$-\tau^2$&$-\tau^2$
&$-1$&$-\tau$&$-1$
&$0$&$0$&$1$
&$1$&$\tau$&$\tau^2$\\
$0$&$-\tau$&$-\tau$
&$-\tau$&$-1$&$-\tau$
&$0$&$-\tau$&$-1$
&$-1$&$0$&$0$
&$1$&$0$&$1$\\
$-\tau$&$-\tau$&$-\tau^2$
&$-1$&$-\tau^2$&$-\tau^2$
&$-\tau$&$-\tau$&$-\tau$
&$0$&$-1$&$0$
&$0$&$1$&$\tau$\\
$-1$&$-\tau$&$-\tau^2$
&$-\tau$&$-\tau^2$&$-2\tau$
&$-\tau$&$-\tau^2$&$-\tau^2$
&$-1$&$-\tau$&$-1$
&$0$&$0$&$1$\\
\hline
\end{tabular}
\end{ex}
\vskip .2cm
The second table in example \ref{e:icos} exhibits certain symmetry properties
which are valid in the general case and which we now address.
\begin{thm}\label{t:tables}
The quantities $\mu_i\cdot \rho_j$ have the following properties.
\begin{description}
\item[\rm (a)] $\mu_i\cdot \rho_{j} = -\mu_{j+n} \cdot \rho_{i}$ for all $i$ and $j$.
\item[\rm (b)] $\mu_i\cdot \rho_j \ge 0$, for $1 \le i \le j \le nh/2$.
\item[\rm (c)] $\mu_{i+t}\cdot \rho_i = 0$, for $1 \le t \le n-1$ and
for all $i$.
\item[\rm (d)] $\mu_{j} \cdot \rho_i \leq 0$ for $1 \le i < j \le nh/2$.
\end{description}
\end{thm}
\paragraph{\it Proof. (a)}
Since $\gamma \mu_i = \mu_i-2\rho_i $ by Theorem
\ref{t:petrieroot},
we compute
\begin{eqnarray*}
-\mu_{j+n}\cdot \rho_i &=& -\gamma \mu_j \cdot (-1/2)(\gamma \mu_i - \mu_i)\\
&=& 1/2(\gamma \mu_j\cdot \gamma \mu_i -\gamma \mu_j\cdot  \mu_i)\\
&=& 1/2(\mu_j\cdot \mu_i -\gamma \mu_j\cdot  \mu_i)\\
&=&  1/2(\mu_j - \gamma \mu_j) \cdot  \mu_i\\
&=&  \rho_j  \cdot \mu_i.
\end{eqnarray*}
(b)  \ \ \ Since $\rho_p\cdot\mu_q = \gamma \rho_p\cdot \gamma \mu_q$
we can assume that $1 \le
i \le n$ and $i < j \le nh/2$.
The result follows from the fact that $\rho_j$ is a positive root and $\mu_i$ is
part of the dual basis to $\{\alpha_1, \dots , \alpha_n\}$.

\noindent
(c) \ \ \ Since $\gamma = R(\rho_{i+t})R(\rho_{i+t-1})\dots
R(\rho_{i})\dots R(\rho_{i+t-n+1})$ (by note \ref{n:consec})
equation~(\ref{e:rhoperpmu}) implies that $\rho_i \cdot \mu_{i+t} = 0$.\\
\noindent
(d) \ \ If $j < i+n$, this follows from part (c).  If $j \ge i+n$,
this property follows from parts (a) and (b).
\vskip .2cm
\begin{cor}\label{c:cones}
If $1 \le i_1 < i_2 < \dots < i_m < k \le nh/2$ , then $\rho_k$
does not lie in the positive cone on $\{\rho_{i_1}, \dots , \rho_{i_m}\}$.
\end{cor}
\emph{Proof:}  If the vector $\vec x$  is expressible as a
non-negative linear combination of $\rho_{i_1}, \dots , \rho_{i_m}$
then $\mu_k\cdot \vec x \le 0$, since $\mu_k\cdot \rho_{i_j} \le 0$,
for $j = 1, \dots , m$.  As $\mu_k \cdot \rho_k = 1$, the vector
$\rho_k$ is not expressible in this manner.
\vskip .2cm
For the construction of the simplicial complex of section $4$ we will
also use the
following result.
\begin{lem}\label{l:flag}
Suppose $\{\sigma_1, \sigma_2, \dots , \sigma_k\}$ is a
consistently ordered
subset of $\{\rho_1, \ldots,\rho_{nh/2}\}$, and
for each $i$ let $\mu(\sigma_i) = -2(\gamma-I)^{-1}\sigma_i$. Then
the following are equivalent.
\begin{description}
\item[\rm (a)] $l(R(\sigma_1)R(\sigma_{2})\dots R(\sigma_k)\gamma) = n-k$.
\item[\rm (b)] $\mu(\sigma_i) \cdot \sigma_j = 0 \mbox{ whenever } i > j$.
\end{description}
\end{lem}
\emph{Proof:}  (a) $\Rightarrow$ (b): \hspace{.2cm}
Let $\delta = R(\sigma_1)R(\sigma_{2})\dots R(\sigma_k)\gamma$
and assume that (a) is true. If $i > j$ then equation~(\ref{e:conj}) gives
\begin{eqnarray*}
 \gamma &=& R(\sigma_k)\dots R(\sigma_i) \dots R(\sigma_j)
 \dots R(\sigma_1) \delta \\
 &=& R(\sigma_i) \dots R(\sigma_j) \dots R(\sigma_1) \delta
      R(\gamma^{-1}\sigma_k) \ldots R(\gamma^{-1}\sigma_{i+1}).
 \end{eqnarray*}
 Since $l(\delta) = n-k$, it follows from equation~(\ref{e:rhoperpmu}) that
 $R(\sigma_j) \leq R(\sigma_i)\gamma$.
 Thus $\sigma_j \in M(R(\sigma_i)\gamma)$ and hence
 $\mu(\sigma_i) \cdot \sigma_j = 0$
 (since $\mu(\sigma_i) \in F(R(\sigma_i)\gamma)$).

\vskip .2cm
(b) $\Rightarrow$ (a):  \ \   Assume that (b) is true.
Then the $k \times k$ matrix $[\mu(\sigma_i)\cdot\sigma_j]$
is upper triangular with ones on the diagonal (corollary~\ref{c:petrieroot} and
part (c) of theorem~\ref{t:tables}).
Thus this matrix is nonsingular
and it follows that  each of $\{\sigma_1,\ldots,\sigma_k\}$
and $\{\mu(\sigma_1),\ldots,\mu(\sigma_k)\}$ is a linearly independent set.

\vskip .2cm
Next we show that $R(\sigma_k)R(\sigma_{k-1})\dots R(\sigma_{i})
\le \gamma$
by reverse induction on
$i$. The case $i=k$ is immediate since $M(\gamma) = {\bf R}^n$. Now
suppose that $i < k$ and that
$R(\sigma_k)R(\sigma_{k-1})\dots R(\sigma_{i+1}) \le \gamma$.
If $ i+1 \le j \le k$ then using equation~(\ref{e:algpo}) and
equation~(\ref{e:posetmap}) we get
\begin{eqnarray*}
R(\sigma_j) &\le& R(\sigma_k)R(\sigma_{k-1})\dots R(\sigma_{i+1})
\le \gamma\\
\Rightarrow R(\sigma_j)\gamma &\ge&
R(\sigma_{i+1})\dots R(\sigma_{k-1})R(\sigma_k)\gamma\\
\Rightarrow F[R(\sigma_j)\gamma] &\subseteq &
F[R(\sigma_{i+1})\dots R(\sigma_{k-1})R(\sigma_k)\gamma]\\
\Rightarrow \mu(\sigma_j) &\in &
F[R(\sigma_{i+1})\dots R(\sigma_{k-1})R(\sigma_k)\gamma]
\end{eqnarray*}
Since the $\mu(\sigma_j)$ are linearly independent and
$F[R(\sigma_{i+1})\dots R(\sigma_{k-1})R(\sigma_k)\gamma]$ has dimension
$k-i$ by the inductive hypothesis, the set
$\{\mu(\sigma_{i+1}), \dots , \mu(\sigma_k)\}$ is a basis for the
subspace
$F[R(\sigma_{i+1})\dots R(\sigma_{k-1})R(\sigma_k)\gamma]$.
Since $\sigma_i$ is orthogonal to each vector in this basis it follows that
\[\sigma_i \in M[R(\sigma_{i+1})\dots R(\sigma_{k-1})R(\sigma_k)\gamma]\]
and hence
$R(\sigma_k) \ldots R(\sigma_{i+1})R(\sigma_{i})\leq \gamma$,
which completes the inductive step.

Finally, since $l(R(\sigma_k) \ldots R(\sigma_{2})R(\sigma_{1})) = k$
(by linear independence of $\sigma_1,\ldots,\sigma_k$)
and $R(\sigma_k) \ldots R(\sigma_{2})R(\sigma_{1}) \leq \gamma$,
we can conclude that the length of
$R(\sigma_1)R(\sigma_{2})\dots R(\sigma_k)\gamma$ is $n-k$.

\section{Definition of the complex~$X(\gamma)$.}\label{s:Xgamma}
\begin{figure}
\centering
\includegraphics[120pt,220pt][490pt,640pt]{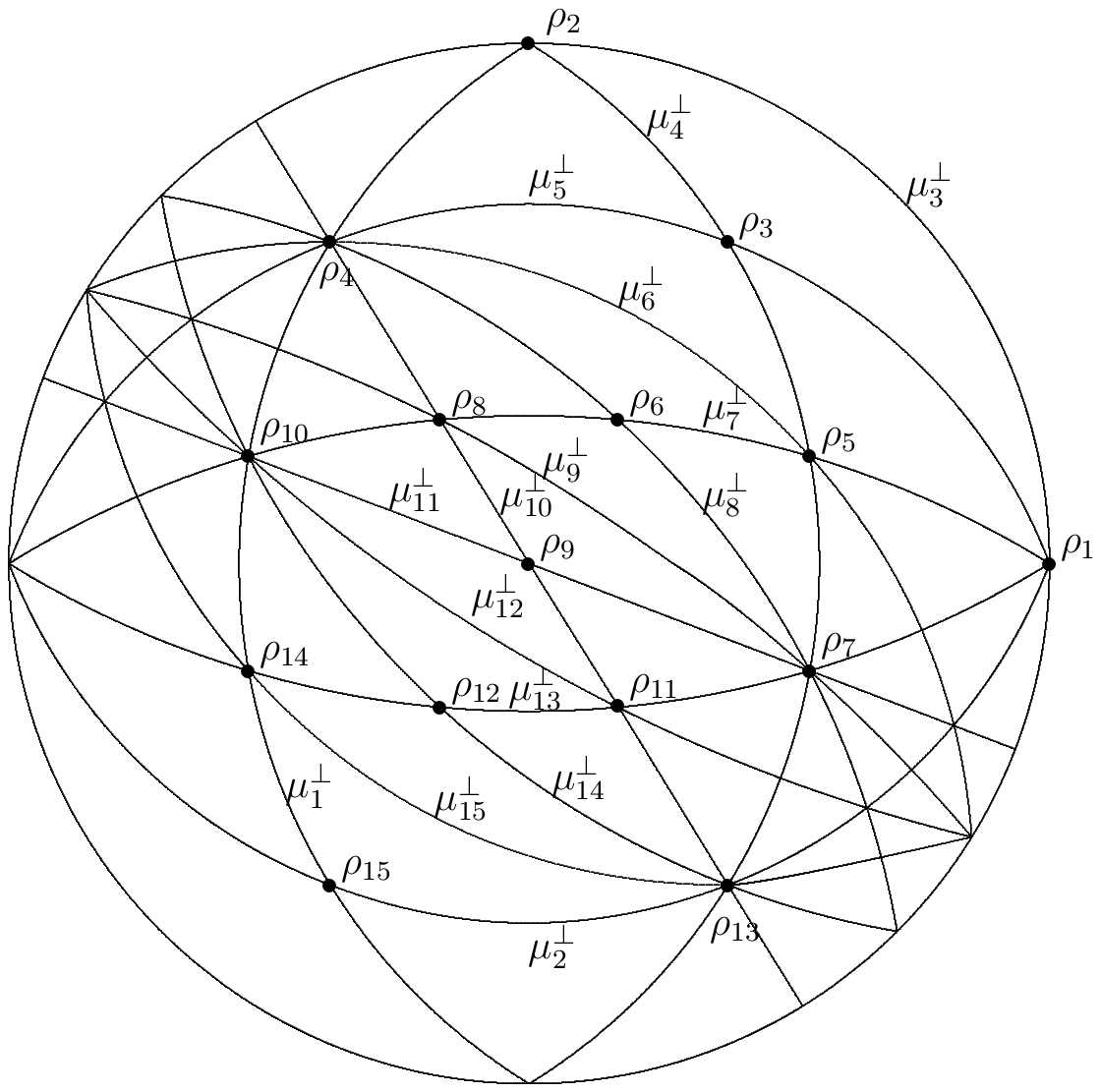}
\caption{}\label{f:icos}
\end{figure}

\begin{figure}
\centering
\includegraphics[120pt,190pt][490pt,630pt]{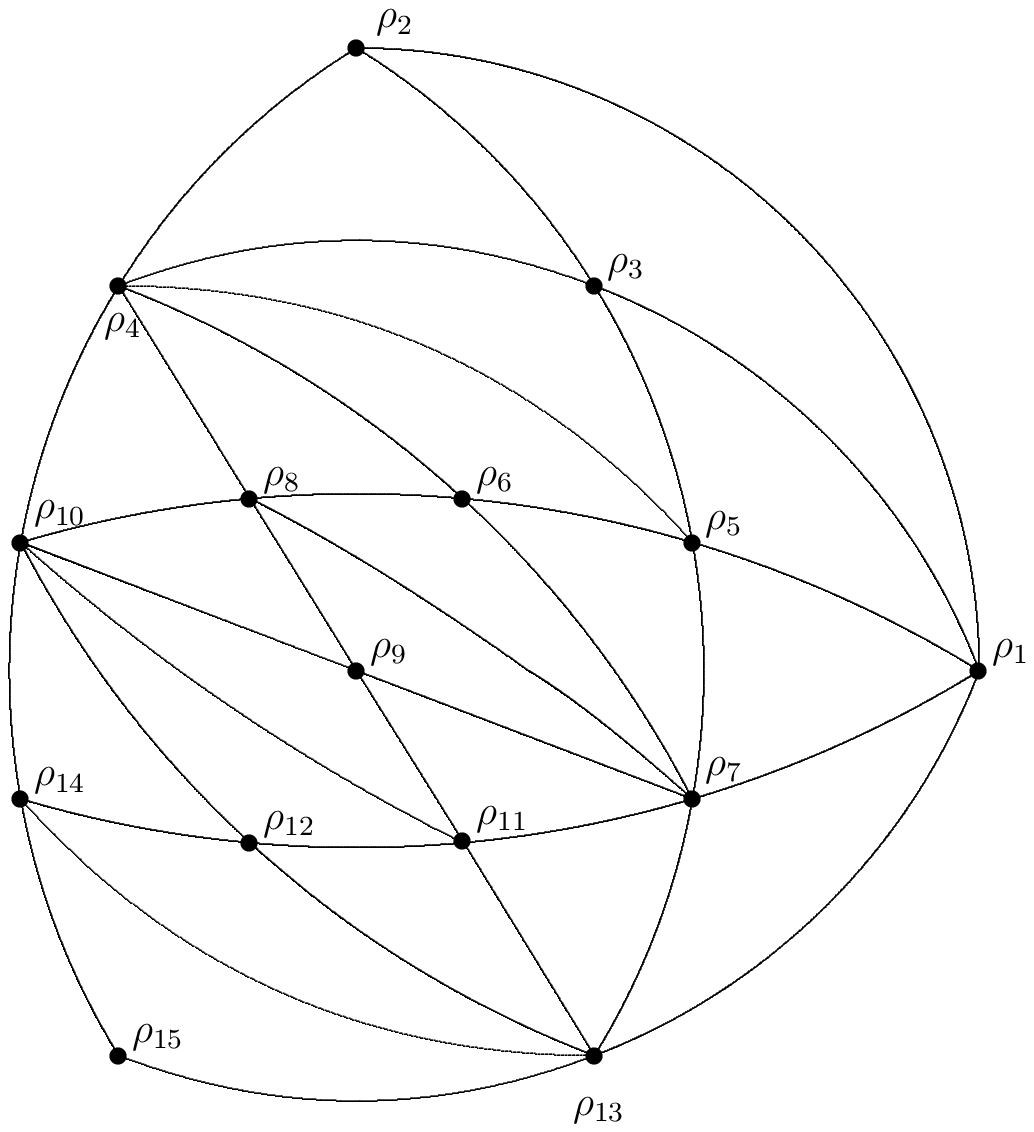}
\caption{}\label{f:icos1}
\end{figure}

Before describing the simplicial complex $X$ which will be used to
prove the lattice property, we give some motivation
for its definition.
By \cite{bwcd}, the map $M$, which associates to an orthogonal
transformation $\alpha$ its moved space $M(\alpha)$, restricts to give
a poset
isomorphism of the interval $[I, \gamma]$ in $W$ onto its
image. Here the set of subspaces of
${\bf R}^n$ is partially
ordered by inclusion and is, in fact, a lattice whose meet operation is
given by subspace intersection. However, the set of subspaces of the form $M(\sigma)$
for $\sigma \le \gamma$ is not closed under intersection.
There are many examples of group elements $\alpha, \beta \le \gamma$ for
which there is no element
$\delta \le \gamma$ in $W$ satisfying $M(\delta) = M(\alpha) \cap
M(\beta)$.  (For example, if $W = A_3 = \Sigma_4$ then the reflections are transpositions
and every $4$-cycle is a Coxeter element.
Consider $\gamma = (1,2,3,4)$ and
the elements $\alpha$ and $\beta$  given by
\[\alpha = (1,3)\gamma = (1,2)(3,4) \mbox{ and } \beta = (2,4)\gamma = (1,4)(2,3). \]
Each of $\alpha$ and $\beta$ has a two-dimensional moved space, these moved spaces
intersect inside the three-dimensional
space $M(\gamma)$ in a one-dimensional subspace, but no transposition precedes
both $\alpha$ and $\beta$.)
To get around this problem we fix a set of
simple roots and replace $M(\alpha)$ by the positive cone on the
positive roots associated
to those reflections which precede $\alpha$. Two surprising things
happen: \emph{(i)}
the intersection of the positive cones associated to two elements
$\alpha$ and $\beta$ of $[I, \gamma]$
is equal to the positive cone associated to some other element
$\delta$ of $[I, \gamma]$ (see proof of theorem~\ref{t:lattice})
and \emph{(ii)} this collection of
 positive cones intersects the unit sphere in a spherical
 simplicial complex which we denote $X(\gamma)$.  We begin this section by
 defining $X(\gamma)$.  Much of our later work is devoted to showing that $X(\gamma)$
 is indeed a simplicial complex.
\vskip .2cm
\begin{defn}\label{d:xgamma}
We define a set, $X = X(\gamma)$, of simplices
by declaring that
\begin{itemize}
\item
the vertex set is $\{\rho_1, \rho_2, \dots , \rho_{nh/2}\}$,
\item
that an edge
joins $\rho_i$ to $\rho_j$ whenever $i < j$ and $R(\rho_i)R(\rho_j)
\le \gamma^{-1}$ and
\item
$\langle \rho_{i_1}, \rho_{i_2}, \dots , \rho_{i_k}\rangle$
forms a $(k-1)$-simplex
if the vertices are distinct and pairwise joined by edges.
\end{itemize}
\end{defn}
\begin{note}\label{n:flag}
It follows from lemma~\ref{l:flag} that a set of $k \leq n$
distinct vertices,
$\{  \rho_{i_1}, \rho_{i_2}, \dots , \rho_{i_k} \}$,
with $1 \le i_1 < i_2 < \dots < i_k \le nh/2$,
determines a $k-1$ simplex of $X$, provided
\[l[R(\rho_{i_1})R(\rho_{i_2})\dots R(\rho_{i_k})\gamma] = n-k.\]
 The largest possible dimension for a simplex of
$X$ is $n-1$  and
at least $nh/2 -n+1$ of these occur by
part (c) of theorem~\ref{t:tables}
and lemma~\ref{l:flag}.
\end{note}
In the case of $I_3$ the complex $X$ is illustrated in
Figures~\ref{f:icos} and \ref{f:icos1}.  The first figure
shows the vertices $\rho_i$ and the moved spaces of the
length two elements in $[I, \gamma]$, while the second
figure shows the simplicial complex $X$.

\vskip .2cm

\section{Subcomplexes}\label{s:subcomplex}
In this and the following sections we will consider a fixed element
$\sigma \le \gamma$.  We let
$P = \{\tau_1, \tau_2, \dots , \tau_t\}$ be the set of positive roots
whose reflections precede $\sigma$.   We assume that $P$ is
ordered consistently with the total order from section $3$ and we let the
corresponding ordered set of vertices of the Petrie polygon
be $\{\mu(\tau_1), \mu(\tau_2), \dots , \mu(\tau_t)\}$.
In this section we will characterise the simple system
associated to $\sigma$ and we will analyse the case where
$l(\sigma) = 2$.  We let $X(\sigma)$ be the collection of simplices of $X(\gamma)$
whose vertices lie in $P$.
\vskip .2cm
Let $k= l(\sigma)$.  Then $P$ is the
set of positive roots for the, possibly reducible, rank $k$
reflection group $\{\alpha \in W \mid M(\alpha) \subseteq
M(\sigma)\}$. Thus $P$ contains a simple system $\Delta =
\{\delta_1, \delta_2, \dots , \delta_k \}$, which we assume is
ordered consistently with the total order from section $3$.
\begin{thm}\label{t:simplesystem}
The ordered elements $\delta_1, \delta_2, \dots , \delta_k $ of
$\Delta$ are determined recursively  by the fact that $\delta_i$
is the last positive root in the subspace $M(\sigma
R(\delta_k)R(\delta_{k-1}) \dots R(\delta_{i+1}))$. In particular,
$\delta_k = \tau_t$.
\end{thm}
\emph{Proof:} Since $\tau_t$ is the last root in $P$ we can apply
Corollary~\ref{c:cones} to $\tau_t$ and the set $P-\{\tau_t\}$ to
deduce that $\tau_t$ cannot be a positive linear combination of
the elements of $P-\{\tau_t\}$. Thus $P-\{\tau_t\}$ cannot contain
a simple system and $\delta_k = \tau_t$.

\vskip .2cm

Next we show that $U_k=M(\sigma R(\delta_k))$ contains
$\{\delta_1,\ldots,\delta_{k-1}\}$. To begin, we show that
 $U_k = M(\sigma) \cap \mu^\perp$, where $\mu$ is the Petrie polygon
 vertex~$\gamma(\mu(\delta_k))$. Note that $U_k =
M(\sigma R(\delta_k))$ has dimension $k-1$. By parts (a) and (b)
of Theorem~\ref{t:tables}, if $1 \le i \le t$ then
\[\mu \cdot \tau_i  = \gamma (\mu(\tau_t)) \cdot \tau_i =
-\tau_t \cdot \mu(\tau_i)\le 0.\] For $i=t$, we obtain $\mu \cdot
\tau_t = -1$. Thus $\tau_t$ is not an element of
$M(\sigma)\cap\mu^\perp$, and this subspace must also have
dimension $k-1$. Since $M(\sigma R(\delta_k))$ is contained in
$M(\sigma)$ and
$$M(\sigma R(\delta_k)) \subseteq M(\gamma
R(\delta_k))=M(R[\gamma(\delta_k)]\gamma)=\mu(\gamma(\delta_k))^\perp
= \gamma(\mu(\delta_k))^\perp$$
(by equation~(\ref{e:knockoff}) and equation~(\ref{e:conj}))
it follows that $U_k \subseteq
M(\sigma) \cap \mu^\perp$, and hence these subspaces are equal.

Now suppose $\tau$ is an arbitrary element of  $P \cap U_k$  and
write
\[\tau = a_1\delta_1 + \dots + a_k \delta_k , \mbox{ with each }a_i \ge 0.\]
Since $\tau \in \mu^\perp$ and $\delta_i \cdot \mu \le 0$ for
$i=1,\ldots,k$, we must have $a_i = 0$ or $\delta_i \cdot \mu = 0$
for each $i$ . Thus $\tau$ is a linear combination of simple roots
in $ M(\sigma)\cap \mu^{\perp}=U_k$. However, the span of $P \cap
U_k$ is all of $U_k$ since $l(\sigma R(\tau_t)) = k-1$. Hence the
$k-1$ simple roots~$\delta_1,\ldots,\delta_{k-1}$ must lie in
$U_k$.

\vskip .2cm

The theorem follows since we can now apply the same arguments as
above to the shorter element $\sigma' = \sigma R(\delta_k)$.

\begin{note}\label{n:delta1}
In an analogous way we can show that
$\delta_i$ is the first root in
$M( R(\delta_1)R(\delta_{2}) \dots R(\delta_{i-1})\sigma)$ and,
in particular, $ \delta_1 = \tau_1$.
\end{note}

\begin{defn}\label{d:fatsimplex}
We will refer to the intersection with $S^{n-1}$ of the positive cone on
$\{\delta_1, \delta_2, \dots , \delta_k \}$ as the fat simplex
associated to $\sigma$.
\end{defn}

\begin{thm}\label{t:Xdih}
If $\sigma \le \gamma$ has length $2$
then $X(\sigma)$ consists of the $t$ roots from $P$ and the
$(t-1)$ $1$-cells given by
$ \langle \tau_i, \tau_{i+1} \rangle$ for $ i = 1,2, \dots , t-1$.
\end{thm}

\emph{Proof:} Suppose $\sigma \le \gamma$ has length two. Since
both $R(\tau_i)\sigma$ and $\sigma R(\tau_i)$ have length one for
each $i$, each subspace~$M(\sigma)\cap \mu(\tau_i)^{\perp} $
contains precisely one positive root $\tau_j$ and each $\tau_j$
lies in precisely one subspace of the form $M(\sigma)\cap
\mu(\tau_i)^{\perp}$. It follows that the $t \times t$ matrix
$A=[\mu(\tau_i)\cdot\tau_j]$ has precisely one zero in each row
and column. Note  that the diagonal entries are all $1$, the
entries  below the diagonal are non-positive, while the entries
above the diagonal are non-negative (by theorem~\ref{t:tables}).
By theorem~\ref{t:simplesystem} and note~\ref{n:delta1},
$\sigma = R(\tau_1)R(\tau_t)$ and
hence  $\tau_t \cdot \mu(\tau_1) = 0$. Let $i_1 > 1$ be the value
for which  $\sigma = R(\tau_{i_1})R(\tau_1)$. Thus
$\mu(\tau_{i_1})\cdot \tau_1 = 0$ and $X(\sigma)$ contains the
edge $\langle \tau_1,\tau_{i_1}\rangle$.  As row $i_1$ contains
only one zero and $i_1 \leq t$, we must have $\mu(\tau_{i_1})\cdot
\tau_t > 0$ (by theorem~\ref{t:tables}). Now, since
$\{\tau_1,\tau_t\}$ is a simple system,
we can write each $\tau_j$ as a non-negative linear combination
of $\tau_1$ and $\tau_t$.  Thus we obtain
$\tau_j \cdot \mu(\tau_{i_1})\ge 0$.
It follows that each entry in row $i_1$ of $A$ is non-negative.  Since every row after
the second has at least two entries below the diagonal and no more than
one of these can be zero we get $i_1 \le 2$.  Since we know that
$i_1 \ne 1$ we deduce that $i_1$ must be equal to $2$. If $t=2$
the proof is complete.
\vskip .2cm
The rest of the proof uses induction.  We assume that $r <t$ and we have
\[\sigma = R(\tau_1)R(\tau_t) = R(\tau_2)R(\tau_1) = \dots = R(\tau_r)R(\tau_{r-1}).\]
We know that $\sigma = R(\tau_{i_r})R(\tau_{r})$ for some $i_r > r $.
Thus
\[\mu(\tau_{i_r})\cdot \tau_r = 0 \mbox{ while } \mu(\tau_{i_r})\cdot \tau_t > 0.\]
If $j \ge r$ then corollary~\ref{c:cones} and the fact that
$M(\sigma)$ is two dimensional implies that $\tau_j$ can be expressed
as a non-negative linear combination
of $\tau_r$ and $\tau_t$.  As in the case $r = 1$, we obtain
$\tau_j \cdot \mu(\tau_{i_r})\ge 0$ for $r \le j \le t$.
It follows that each of the last $t-r+1$ entries in row $i_r$ of $A$ is non-negative.
Since row $j$ has $j-1$ entries below the diagonal and no more than
one of these can be zero we get $i_r \le  r+1$.  Since we know that
$i_r > r$ we deduce that $i_r$ must be equal to $r+1$.
\vskip .2cm

In the case $l(\sigma) = 2$, the action of the simple reflections
on $P$ can now be deduced.

\begin{lem}\label{l:simpleaction}
If $l(\sigma) = 2$ and
$P = \{\tau_1, \dots , \tau_t\}$ is the ordered set of positive roots
for reflections preceding $\sigma$ then
\[R(\tau_1)(\tau_i) = \left\{
\begin{array}{ll}
\tau_{t-i+2} &\mbox{ for } 2 \le i \le t\\
-\tau_1 &\mbox{ for }i = 1,
\end{array}
\right. \]
\[\mbox{ and } R(\tau_t)(\tau_i) = \left\{
\begin{array}{ll}
\tau_{t-i} &\mbox{ for } 1 \le i \le t-1\\
-\tau_t &\mbox{ for }i = t.
\end{array}
\right.\]
\end{lem}

\emph{Proof:}
We will consider the case of $\tau_1$.  The case of
$\tau_t$ is similar.
From the proof of theorem~\ref{t:Xdih} it follows  that the only
expressions for $\sigma$ as a product of two reflections are
\[\sigma = R(\tau_1)R(\tau_t) = R(\tau_2)R(\tau_1) =
R(\tau_3)R(\tau_2) = \dots = R(\tau_t)R(\tau_{t-1}).\] Since
$R(\tau_1)$ permutes $P-\{\tau_1\}$ (lemma A in I.4.3 of
\cite{Kane}) and
\[ \sigma = R(\tau_1)R(\tau_t) = (R(\tau_1)R(\tau_t)R(\tau_1))R(\tau_1)
= R(R(\tau_1)\tau_t) R(\tau_1)\]
we deduce that $R(\tau_1)\tau_t=\tau_2$. In general, suppose that
 $1 < i+1 \leq t$ and that $R(\tau_1)\tau_{i+1} = \tau_{t-(i+1)+2}$.
 Conjugate $R(\tau_t)R(\tau_1) = \sigma^{-1}
=  R(\tau_i)R(\tau_{i+1})$ by $R(\tau_1)$ to obtain
\[ R(\tau_1)R(\tau_t) = \sigma =  R(R(\tau_1)\tau_i)R(R(\tau_1)\tau_{i+1})
= R(R(\tau_1)\tau_{i})R(\tau_{t-i+1}).\]
Comparing with the list of factorisations for $\sigma$, we deduce that
$R(\tau_1)\tau_{i} = \tau_{(t-i+1)+1}=\tau_{t-i+2}$. The result follows by reverse
induction on $i$ starting at $i = t$.

\begin{lem}\label{l:rootangle}
Let $\rho_i$ and $\rho_j$ be  distinct positive roots
for which  $R(\rho_i)R(\rho_j) \le \gamma$.\\
(a) If $i < j$ then $\rho_i \cdot \rho_j \le 0$.\\
(b) If $i > j$ then $\rho_i \cdot \rho_j \ge 0$.
\end{lem}

\emph{Proof:}
Let $P = \{\tau_1, \dots , \tau_t \}$ be the ordered set
of positive roots whose reflections precede $\sigma = R(\rho_i)R(\rho_j)$.
 By theorem~\ref{t:Xdih}, $\{\tau_1, \tau_t\}$ is the corresponding simple system
 and the  only expressions for $\sigma$ as a product of
two reflections are
\[\sigma = R(\tau_1)R(\tau_t) = R(\tau_2)R(\tau_1) =
R(\tau_3)R(\tau_2) = \dots = R(\tau_t)R(\tau_{t-1}).\]
Since $\sigma = R(\tau_1)R(\tau_t)$ is the only case
in which the roots appear in increasing order, if $i<j$ then we must have
 $\rho_i=\tau_1$ and $\rho_j=\tau_t$. Part (a) now follows from
$\rho_i\cdot \rho_j = \tau_1 \cdot \tau_t \le 0$ (since $\{\tau_1,\tau_t\}$ is
a simple system).

If $i>j$ then $\rho_i = \tau_{r+1}$ and $\rho_{j} = \tau_{r}$ for
some $r \geq 1$. Thus we need to verify that
$\tau_i\cdot\tau_{i+1} \ge 0$ for $i=1,\ldots,t-1$. Using
lemma~\ref{l:simpleaction} and the facts that reflection is an
isometry and $\{\tau_1,\tau_t\}$ is a simple system, we obtain
\[\tau_1 \cdot \tau_2 = R(\tau_1)\tau_1 \cdot R(\tau_1)\tau_2 =
-\tau_1 \cdot \tau_t \ge 0.\]
For the same reasons, if $i \geq 2$ then
\[\tau_i \cdot \tau_{i+1} = R(\tau_t)R(\tau_1)\tau_i \cdot R(\tau_t)R(\tau_1)\tau_{i+1} =
R(\tau_t)\tau_{t-i+2} \cdot R(\tau_t)\tau_{t-i+1} = \tau_{i-2} \cdot \tau_{i-1} \]
and the result follows by induction on $i$.

\begin{note}\label{n:obtuse}
It is immediate from lemma~\ref{l:rootangle} that the edges
of $X$ can
now be characterised by the following geometric criterion.
There is an edge joining the distinct vertices $\rho_i$ and $\rho_j$
if and only if the vectors subtend a
non-obtuse angle and one of
$R(\rho_i)R(\rho_j)$ or $R(\rho_j)R(\rho_i)$ precedes
$\gamma$.
\end{note}

\section{Walls of  fat simplices}\label{s:walls}
In this section we continue to investigate a fixed element $\sigma$ of
length $k$
which precedes $\gamma$ in $W$. Using the same notation as in the preceding section,
we find a dual basis to $\{\delta_1, \dots , \delta_k\}$
(corollary~\ref{c:dual}).
We also determine the first top-dimensional
simplex of $X(\sigma)$ in the lexicographic order.

\begin{defn}\label{d:epsilon}
For each $ i =1,\ldots k$ we define $\epsilon_i$ to be the root given by
\[\epsilon_i = R(\delta_1)\dots R(\delta_{i-1})\delta_i.\]
\end{defn}

\begin{prop}\label{p:+epsilon}
For each $i=1,\ldots,k$, the vector $\epsilon_i$ is a positive root.
Moreover
\[ \sigma = R(\epsilon_k)R(\epsilon_{k-1})\ldots R(\epsilon_1)\]
and for each $i=1,\ldots,k$
\[\sigma = R(\epsilon_i)R(\delta_1) \dots R(\delta_{i-1})R(\delta_{i+1})
\dots R(\delta_k).\]
\end{prop}

\emph{Proof:} The first statement is a special case of lemma D in section I.4.3
of \cite{Kane}. The remaining statements follow from
\begin{eqnarray*} R(\epsilon_i) &=& R(R(\delta_1)\dots R(\delta_{i-1})\delta_i)\\
&=&[R(\delta_1)R(\delta_2)\dots R(\delta_{i-1})]R(\delta_i)
[R(\delta_{i-1})\dots R(\delta_{2})R(\delta_1)]
\end{eqnarray*}
since each reflection has order $2$.

\bigskip

\begin{corollary}\label{c:walls}
The  walls of the spherical simplex on the vertices
$\{\delta_1, \dots, \delta_k\}$ in the
subspace $M(\sigma)$ are given by intersecting the planes
$\mu^{\perp}(\epsilon_i)$ with the unit sphere in $M(\sigma)$.
\end{corollary}

\emph{Proof:} The second identity in proposition~\ref{p:+epsilon}
implies that the set
$\{\delta_1, \dots, \delta_{i-1}, \delta_{i+1}, \dots , \delta_k\}$
is contained in the subspace $\mu^{\perp}(\epsilon_i)$, for $i = 1,\ldots,k$.

\begin{note}\label{f:fatsimplex}
We will see in the next section that this simplex is the geometric realisation
of $X(\sigma)$.
\end{note}

\begin{defn}\label{d:muprime}
For each $i=1,\ldots,t$, let
$\mu'(\tau_i)$ be
the orthogonal projection
 of $\mu(\tau_i)$ into $M(\sigma)$.
\end{defn}

\begin{prop}\label{p:muprime}
For each $i,j \in \{1, \dots , t\}$, we have
$\mu'(\tau_i)\cdot \tau_j = \mu(\tau_i)\cdot \tau_j$
and $\mu'(\tau_i) \in F(R(\tau_i)\sigma)$.
In particular,
\[\sigma(\mu'(\tau_i)) = R(\tau_i)\mu'(\tau_i) = \mu'(\tau_i)-2 \tau_i.\]
\end{prop}

\emph{Proof:}
Write $\mu(\tau_i) = \mu'(\tau_i) + \vec  y$,
where $\mu'(\tau_i) \in M(\sigma)$ and $ \vec y \in
[M(\sigma)]^{\perp} = F(\sigma)$.
Then $\mu(\tau_i)\cdot \tau_j = \mu'(\tau_i)\cdot \tau_j$
for any $1 \le j \le t$.
Now $\mu(\tau_i) \in F(R(\tau_i)\gamma) \subset F(R(\tau_i)\sigma)$,
because $R(\tau_i)\sigma \le R(\tau_i)\gamma$  by equation~(\ref{e:knockoff}).
But $\vec y $
is also an element of $F(R(\tau_i)\sigma)$, since
$F(R(\tau_i)\sigma)$ contains $F(\sigma)$ by equation~(\ref{e:posetmap}).
It follows that
$\mu'(\tau_i)$ is an element of $F(R(\tau_i)\sigma)$.
The final claim of the proposition follows from this.

\begin{prop}\label{p:dual}  For $1 \le i \le k$, we have
$\mu(\epsilon_i)\cdot \delta_i = 1$.
\end{prop}

\emph{Proof:}
Fix $i$ and let
$\sigma' = R(\delta_1)\dots R(\delta_i)$ (which precedes $\sigma$ and hence
$\gamma$). Write $\mu(\epsilon_i) =
\mu''(\epsilon_i)+ \vec z$ where $\mu''(\epsilon_i) \in M(\sigma')$ and
$\vec z \in M(\sigma')^{\perp} = F(\sigma')$. Applying
proposition~\ref{p:muprime} to $\sigma'$ yields,
\begin{eqnarray*}
\mu(\epsilon_i)\cdot \delta_i &=&
\sigma'[\mu(\epsilon_i)]\cdot \sigma'[\delta_i]\\
&=& \sigma'[\mu''(\epsilon_i)+\vec z]\cdot \sigma'[\delta_i]\\
&=& [\mu''(\epsilon_i)-2\epsilon_i +\vec z]\cdot
R(\delta_1)\dots R(\delta_{i-1})R(\delta_i)[\delta_i] \\
&=& [\mu(\epsilon_i)-2\epsilon_i ]\cdot
R(\delta_1)\dots R(\delta_{i-1})[-\delta_i]\\
&=& [\mu(\epsilon_i)-2\epsilon_i ]\cdot [-\epsilon_i]\\
&=& 1 \mbox{ \ \ as required.}
\end{eqnarray*}

The following corollary is immediate from corollary~\ref{c:walls} and
proposition~\ref{p:dual}

\begin{cor}\label{c:dual}
The dual basis to $\{\delta_1,\ldots,\delta_k\}$ is
$\{\mu'(\epsilon_1),\ldots,\mu'(\epsilon_k)\}$.
\end{cor}

\vskip .2cm
The order $\epsilon_1,\ldots,\epsilon_k$ need not be
consistent with the global ordering. However the induced order
still determines a factorisation of $\sigma$
(proposition~\ref{p:commroot} below). First we make a definition.
\begin{defn}\label{d:theta}
Let $\theta_1, \dots , \theta_k$ be the reordering of
 $\epsilon_1,\ldots,\epsilon_k$
which is consistent with the global order.
\end{defn}

\begin{lem}\label{l:commroot}
If $i < j$ but $\epsilon_i > \epsilon_j$ then
$R(\epsilon_j)R(\epsilon_i)= R(\epsilon_i)R(\epsilon_j)$.
\end{lem}

\emph{Proof:} Assume that $i<j$ and that $\epsilon_i > \epsilon_j$.
We need to show that $\epsilon_i \cdot \epsilon_j = 0$.
Since $\{\delta_1, \dots  , \delta_k \}$
is a simple system, each of the dot products $\delta_i\cdot \delta_j$
is non-positive for $i \neq j$. It follows that $R(\delta_{j-1})\delta_j$
is a non-negative linear combination of $\delta_{j-1}$ and $\delta_j$ and,
by induction, that $R(\delta_{i+1})\dots R(\delta_{j-1})\delta_j$ is a
non-negative linear combination of $\delta_{i+1},\delta_{i+2},\dots,\delta_j$.
Hence, $\delta_i \cdot R(\delta_{i+1})\dots R(\delta_{j-1})\delta_j \le 0$.
Now, since $i<j$, we can compute that
\begin{eqnarray*}
\epsilon_i \cdot \epsilon_j &=&
R(\delta_1)\dots R(\delta_{i-1})\delta_i \cdot
R(\delta_1)\dots R(\delta_{j-1})\delta_j\\
&=& \delta_i \cdot R(\delta_i)\dots R(\delta_{j-1})\delta_j\\
&=& R(\delta_i)\delta_i \cdot R(\delta_{i+1})\dots R(\delta_{j-1})\delta_j\\
&=& -\delta_i \cdot R(\delta_{i+1})\dots R(\delta_{j-1})\delta_j\\
&\ge & 0.
\end{eqnarray*}
However, since $i < j$,  the identity
$\sigma = R(\epsilon_k) \dots R(\epsilon_j) \dots R(\epsilon_i) \dots R(\epsilon_1)$
implies that
$R(\epsilon_j)R(\epsilon_i)$ precedes $\sigma$ and hence $\gamma$. Now
part (a) of
lemma~\ref{l:rootangle} implies that $\epsilon_i\cdot \epsilon_j \leq 0$
(because $\epsilon_i >\epsilon_j$).
We conclude that
$\epsilon_i \cdot \epsilon_j = 0$, as required.

\begin{prop}\label{p:commroot}
The elements $\theta_1,\ldots,\theta_k$ satisfy the identity
\[\sigma =R(\theta_k)\dots R(\theta_1).\]
\end{prop}

\emph{Proof:} We know that $\sigma =R(\epsilon_k)\dots R(\epsilon_1)$.
If $\epsilon_i=\theta_i$ for each $i$, there is nothing to prove.
Otherwise, repeated application of lemma~\ref{l:commroot} yields the required result.

\begin{cor}\label{c:firstsimplex}
The simplex $\langle \theta_1, \dots, \theta_k \rangle$ is
the first top dimensional simplex of $X(\sigma)$ in the lexicographic order.
\end{cor}

\emph{Proof:  } Since $\theta_1<\theta_2<\cdots < \theta_k$ and
$\sigma = R(\theta_k)\dots R(\theta_1)$,  the $(k-1)$-simplex
$\langle \theta_1, \dots, \theta_k \rangle$ is in $X(\sigma)$.
As $\{\delta_1,\ldots,\delta_k\}$ is the simple system corresponding to $P$,
$\{\theta_1, \dots, \theta_k\}$ is a rearrangement of
$\{\epsilon_1, \dots , \epsilon_k\}$ and
\[
\mu(\epsilon_a)\cdot \delta_b = \left\{
\begin{array}{ll}
0 &\mbox{ for } a \ne b\\
1 &\mbox{ for }a = b,
\end{array}
\right.
\]
it follows that  $\mu(\theta_i)\cdot \tau_j \ge 0$
for $1 \le i \le k$ and $1 \le j \le t$. However,  theorem~\ref{t:tables}
implies that whenever $\tau_j < \theta_i$,
the dot product $\mu(\theta_i)\cdot \tau_j$ is non-positive, and hence it
must be zero. Therefore
\[ \{\tau_j \mid \tau_j < \theta_{i}\} \subseteq
M(\sigma)\cap \mu(\theta_i)^\perp\cap \mu(\theta_{i+1})^\perp \cap
\dots \cap \mu(\theta_k)^\perp \]
which is an $(i-1)$-dimensional subspace of ${\bf R}^n$ since the $\theta$'s are
linearly independent.
Now if $\langle \tau_{i_1}, \dots, \tau_{i_k} \rangle$ is a $(k-1)$-simplex
of $X(\sigma)$ with $\tau_{i_1} < \tau_{i_2} < \dots < \tau_{i_k}$,
then $\langle \tau_{i_1}, \dots, \tau_{i_j} \rangle$ is a
$(j-1)$-simplex for each $j \le k$ and this forces $\tau_{i_j} \ge \theta_j$.
\vskip .2cm
We finish this section with two examples.  The first shows that
the first $k$ roots of $P$ may not span a top dimensional simplex.
The second example illustrates that even in the case $\sigma = \gamma$ the
$\epsilon$'s may not be ordered consistently with the global order.

\begin{figure}
\centering
\includegraphics[90pt,270pt][500pt,610pt]{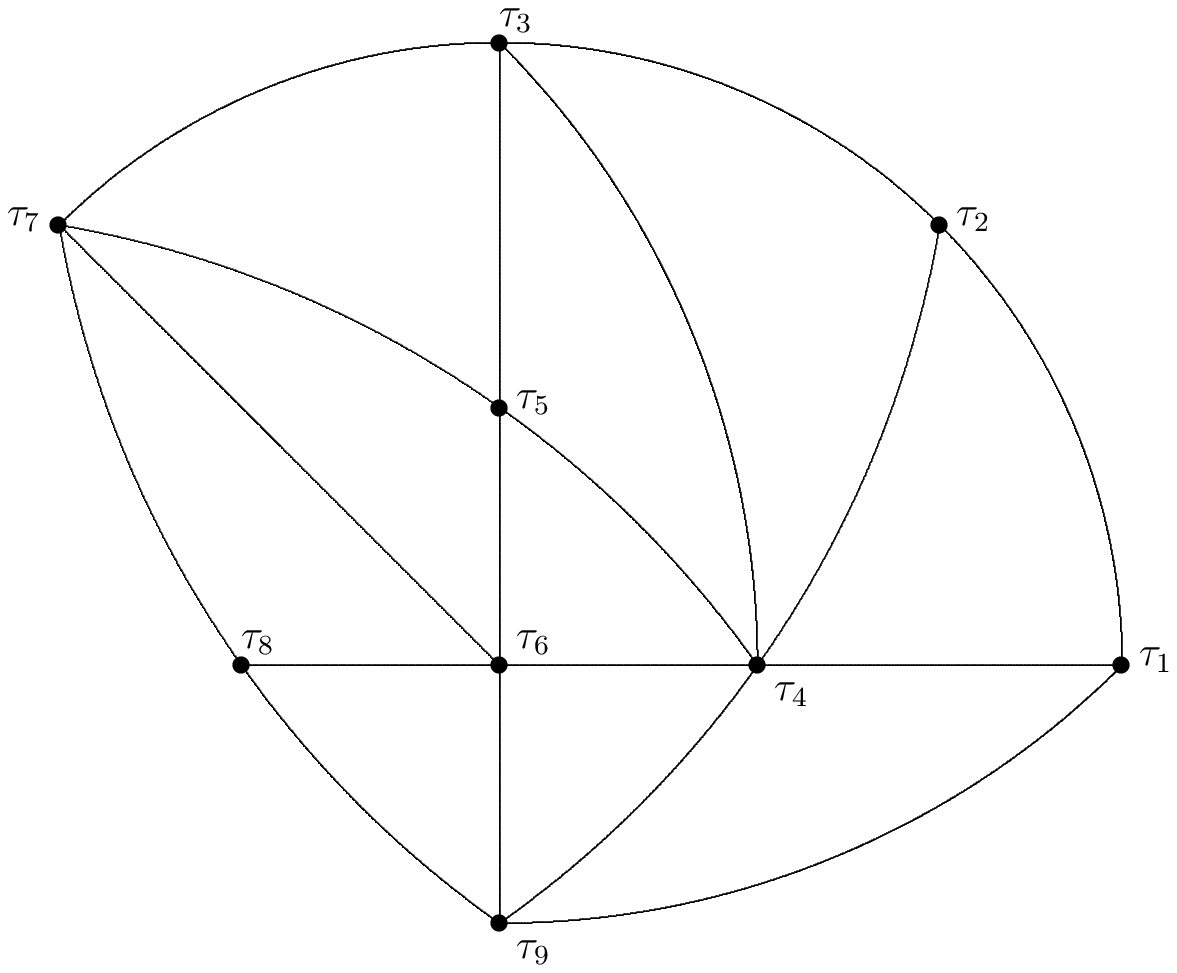}
\caption{}\label{f:cube}
\end{figure}

\begin{ex}\label{e:cube}
Consider the symmetry group of the $4$-dimensional cube.
One simple system of unit vectors for this group is
\[\alpha_1 = (1,0,0,0), \alpha_2 = (\sqrt{2}/2)(0,1, 0,-1),\]
\[\alpha_3 = (\sqrt{2}/2)(-1,0, 0,1), \alpha_4 = (\sqrt{2}/2)(0,-1,1,0).\]
The element  $\gamma = [1,2,3,4]$ (in
the notation of \cite{bwcd}) is one of the Coxeter
elements determined by this simple system,
where
\[[1,2,3,4](x,y,z,w) = (-w,x,y,z).\]
The dual basis to $\{\alpha_1,\alpha_2,\alpha_3,\alpha_4\}$ is
\[\beta_1 = (1,1,1,1), \beta_2 = (\sqrt{2})(0,1, 1,0),\]
\[\beta_3 = (\sqrt{2})(0, 1,1,1), \beta_4 =  (\sqrt{2})(0,0,1,0).\]
The $16$ positive roots and  the $16$ corresponding Petrie polygon vertices
can be computed according to section~\ref{s:root-mu}.
The element $\sigma = [1,2,3]$ precedes $\gamma$, has length three
and its nine positive roots are tabulated
below.
\vskip .2cm
\begin{tabular}{|l|l|l|l|l|}
\hline
$i$&$\tau_i$&$\rho$ subscript&$\mu(\tau_i)$\\
\hline
$1$&$(1,0, 0,0)$&$1$&$(1,1,1,1)$\\
$2$&$(\sqrt{2}/2)(1,1, 0,0)$&$3$&$(\sqrt{2})(0,1,1,1)$\\
$3$&$(0,1,0,0)$&$5$&$(-1,1,1,1)$\\
$4$&$(\sqrt{2}/2)(1,0, 1,0)$&$6$&$(\sqrt{2})(0,0,1,1)$\\
$5$&$(\sqrt{2}/2)(0,1, 1,0)$&$7$&$(\sqrt{2})(-1,0,1,1)$\\
$6$&$(0,0,1,0)$&$9$&$(-1,-1,1,1)$\\
$7$&$(\sqrt{2}/2)(-1,1, 0,0)$&$12$&$(\sqrt{2})(-1,0,0,0)$\\
$8$&$(\sqrt{2}/2)(-1,0, 1,0)$&$14$&$(\sqrt{2})(-1,-1,0,0)$\\
$9$&$(\sqrt{2}/2)(0,-1, 1,0)$&$15$&$(\sqrt{2})(0,-1,0,0)$\\
\hline
\end{tabular}
\vskip .2cm
We use theorem~\ref{t:simplesystem} to determine the
simple system~$\{\delta_1,\delta_2,\delta_3\}$.
First $\delta_3 = \tau_9$ and hence $\sigma R(\delta_3) = [1,2]$.  It follows that
$\delta_2 = \tau_7$ and, since
$\sigma R(\delta_3)R(\delta_2) = [1]$, that $ \delta_1 = \tau_1$.
Thus
\[\Delta = \{(1,0,0,0),(\sqrt{2}/2)(-1,1, 0,0) ,
(\sqrt{2}/2)(0,-1, 1,0) \}.\]
Next we calculate
\[\epsilon_1 = \delta_1 = \tau_1,\,\,\,
\epsilon_2 = R(\delta_1)\delta_2 = \tau_2, \mbox{ \ and \ }
\epsilon_3 = R(\delta_1)R(\delta_2)\delta_3 = \tau_4.\]
The ordering $\epsilon_1,\epsilon_2,\epsilon_3$ is already
subordinate to the global ordering.  Thus $\theta_i = \epsilon_i$ in this case.
Figure~\ref{f:cube} shows
the subcomplex $X(\sigma)$ for this example.
\end{ex}
\vskip .2cm

\begin{ex}\label{e:tetra}
Consider the group~$A_3 = \Sigma_4$.  We have seen that the four-cycle
$\gamma = (1,2,3,4)$ is a Coxeter element.  Choose the simple system
\[\alpha_1 = (\sqrt{2}/2)(1,0,-1,0), \alpha_2 = (\sqrt{2}/2)(-1,1,0,0),
\alpha_3 = (\sqrt{2}/2)(0,0,1,-1),\]
corresponding
to the factorisation $\gamma = (1,3)(1,2)(3,4)$.  For notational convenience
we will identify roots with the transpositions they determine.  The global order
is the following.
\[(1,3), (2,3), (1,4), (2,4), (3,4), (1,2).\]
We use  theorem~\ref{t:simplesystem} to find the simple system.  First
$\delta_3 = (1,2)$ since this is the last root.  Since $(1,2,3,4)(1,2) = (1,3,4)$,
$\delta_2 = (3,4)$.  Finally, $\delta_1 = (1,3)$.
Thus $\epsilon_1 = (1,3) = \tau_1$, $\epsilon_2 = (1,4) = \tau_3$ and
$\epsilon_3 = (2,3) = \tau_2$, making the ordering on the
$\epsilon$'s inconsistent with the global order.
\end{ex}

\section{Characterisation of $|X(\sigma)|$ and
proof of the lattice property.} \label{s:lattice}
In this section we show that if $\sigma \le \gamma$
then $X(\sigma)$ is a simplicial complex.
We characterise the geometric realisation~$|X(\sigma)|$ of the
subcomplex $X(\sigma)$ and use this characterisation to prove that
the interval $[I, \gamma]$ in $W$ is a lattice.  Throughout this section we continue to
use the notation of the earlier sections.

\vskip .2cm
We begin with a technical result (proposition~\ref{p:separate})
about the
separation properties of the hyperplanes $\{\mu(\tau_i)^{\perp}\}$
which is used in the proof of theorem~\ref{t:X = Z}.
It depends on proposition~\ref{p:sigmaaction} which is stated for
convenient reference and which concerns
the action of $\sigma^{-1}$ on $P$.
Using the factorisation $\sigma^{-1} = R(\delta_k)\dots R(\delta_1)$,
this is a special case of theorem B in section I.4.3 of~\cite{Kane}.

\begin{prop}\label{p:sigmaaction}
If $\tau_s \in P$ then
$\sigma^{-1}(\tau_s) \in -P$ if and only if
$\tau_s = \epsilon_i$ for some $i$ with $1 \le i \le k$.
\end{prop}

\begin{prop}\label{p:separate}
Let $\tau_a$ and $\tau_s$ be elements of $P$ such that
$\tau_a \ge \theta_k$,  $\tau_s < \tau_a$,
$\tau_s \notin \{\theta_1,\theta_2,\dots,\theta_{k-1},\theta_k\}$ and
$\tau_a \in \mu(\tau_s)^{\perp}$.   Then we can find two roots
$\tau_b$ and $\tau_c$ in $P$ with $\tau_b, \tau_c < \tau_a$, which both
lie on the hyperplane $\mu(\tau_a)^{\perp}$, but which are separated
by the hyperplane $\mu(\tau_s)^{\perp}$.
\end{prop}

\emph{Proof:}
Since $\tau_a \cdot \mu(\tau_s) = 0$, $\tau_a \in M(R(\tau_s)\gamma)$ and the element
$R(\tau_s)R(\tau_a)$ precedes $\gamma$ by equation~(\ref{e:calc}).

\vskip .2cm
\textbf{Construction of $\tau_b$:  }
Let $Q=  \{ \tau_i \in P \mid R(\tau_i) \le R(\tau_s)R(\tau_a)\}$.
Since $\tau_s < \tau_a$, the proof of lemma~\ref{l:rootangle} implies that (i)
the set $\{\tau_s, \tau_a\}$ is a simple system  which spans
 $Q$ and (ii)  $\tau_s$ is the first root
and $\tau_a$ the last root in the induced ordering of $Q$.  By lemma~\ref{l:simpleaction},
the element $\tau_b$ given by $\tau_b = R(\tau_a)\tau_s$ is in $Q$
 and $s < b < a$.  Since
\[ R(\tau_a)R(\tau_b) = R(\tau_a)(R(\tau_a)R(\tau_s)R(\tau_a)) =
R(\tau_s)R(\tau_a) \le \gamma,\]
we deduce  that $\tau_b \cdot \mu(\tau_a) = 0$.
Finally  $\tau_b\cdot \mu(\tau_s) >0$ because
\[\mu(\tau_s)\cdot \tau_b = \mu(\tau_s)\cdot R(\tau_a)\tau_s
= \mu(\tau_s)\cdot \left\{\tau_s -2(\tau_a \cdot \tau_s)\tau_a \right\}= 1\]
since $\mu(\tau_s)\cdot\tau_s = 1$ and $\mu(\tau_s)\cdot \tau_a = 0$.
\vskip .2cm
\textbf{Construction of $\tau_c$:  }
As $ \tau_s$ is not an element of
$\{\theta_1, \dots, \theta_k\}= \{\epsilon_1,\ldots,\epsilon_k\}$,
the root~$\tau_c$ defined by
$\tau_c = \sigma^{-1}\tau_s$ is an element of $P$ (by
proposition~\ref{p:sigmaaction}).
Since  $R(\tau_s)R(\tau_a) \le \gamma$ and $M[R(\tau_s)R(\tau_a)]
\subset M(\sigma)$ we deduce by equation~(\ref{e:induce}) that $R(\tau_s)R(\tau_a) \le \sigma$.
This gives $R(\tau_a) \le R(\tau_s)\sigma = \sigma R(\tau_c)$
by equation~(\ref{e:conj})
and
hence $R(\tau_a)R(\tau_c) \le  \sigma \le \gamma$
by equation~(\ref{e:calc}).  It follows that
$\tau_c \cdot \mu(\tau_a) = 0$.  Now $\tau_c\cdot \mu(\tau_s) < 0$ since
\begin{eqnarray*}
\mu(\tau_s) \cdot \tau_c &=& \mu'(\tau_s) \cdot \tau_c \\
&=& \mu'(\tau_s) \cdot \sigma^{-1}(\tau_s)\\
&=& \sigma(\mu'(\tau_s)) \cdot \tau_s\\
&=& \{\mu'(\tau_s)-2\tau_s\} \cdot \tau_s \mbox{ \ by proposition~\ref{p:muprime}}\\
&=& 1-2  = -1.
\end{eqnarray*}
Since $\mu(\tau_s) \cdot \tau_c < 0$, theorem~\ref{t:tables} implies that
$\tau_c < \tau_s$ and hence $\tau_c < \tau_a$, as required.

\begin{defn}\label{d:filter}
If $\tau_i$ is a root in $P$ and
$\rho$ is any positive root,  we define
\begin{itemize}
\item $\mu(\tau_i)^+ = \{x\in {\bf R^{n}} \mid x \cdot \mu(\tau_i) \ge 0\}$ (a positive
halfspace),
\item $\mu(\tau_i)^- = \{x\in {\bf R^{n}} \mid x \cdot \mu(\tau_i) \le 0\}$ (a negative
halfspace),
\item $X(\sigma, \rho) = $ the set of simplices of $X$ whose vertices both lie in
$M(\sigma)$ and precede $\rho$ in the total order,
\item $c\left [X(\sigma , \rho) \right ] = $ the positive cone on  the set
$|X(\sigma, \rho)|$.
\item $Y(\sigma, \rho)  = $ the positive cone on  those roots which both lie in
$M(\sigma)$ and precede $\rho$ in the total order,
\item $Z(\sigma, \tau_i) = M(\sigma) \cap \mu(\theta_1)^+ \cap \dots
\cap \mu(\theta_k)^+ \cap \mu(\tau_{i+1})^- \cap \dots \cap \mu(\tau_t)^-$,
for $\tau_i \ge \theta_k$.
\end{itemize}
\end{defn}

\begin{thm}\label{t:Xsimplicial}
The set $X(\sigma, \tau_i)$ is a simplicial complex for each $i=1,\dots,t$.
\end{thm}

\emph{Proof:} We use induction on $i$. First note that
$X(\sigma, \tau_1) = \{ \langle \tau_1 \rangle\}$, a zero-dimensional
simplicial complex.

Assume now that $i\ge 1$ and that
$X(\sigma, \tau_i)$ is a simplicial complex. By definition,
if $\tau_j \in X(\sigma, \tau_i)$ then $\langle \tau_j ,\tau_{i+1}\rangle
\in X(\sigma, \tau_{i+1})$ if and only if $\tau_j \cdot \mu(\tau_{i+1}) = 0$.
It follows from theorem~\ref{t:tables} that the only vertices of
$X(\sigma, \tau_i)$  that are not joined to $\tau_{i+1}$ by an edge in
$X(\sigma, \tau_{i+1})$, lie in the interior of the halfspace
$\mu(\tau_{i+1})^+$.
Hence, $\mu(\tau_{i+1})^\perp \cap |X(\sigma, \tau_i)|$ is a simplicial complex.
Now each simplex in $X(\sigma,\tau_{i+1}) \setminus X(\sigma, \tau_i)$
is of the form $\langle \tau_{a_1},\ldots,\tau_{a_b},\tau_{i+1} \rangle$
where $\tau_{a_1} < \cdots < \tau_{a_b} < \tau_{i+1}$ and where
$\tau_{a_c}\in \mu(\tau_{i+1})^\perp$ for $c=1,\dots,b$. Thus the
simplex
$\langle \tau_{a_1},\ldots,\tau_{a_b} \rangle$ of $X(\sigma,\tau_i)$
is contained in $\mu(\tau_{i+1})^\perp$. Using this, it is straightforward
to verify that the intersection of any two simplices in $X(\sigma,\tau_{i+1})$
is itself a simplex in $X(\sigma,\tau_{i+1})$.

\begin{cor}\label{c:Xsimplicial}
For each $\sigma \le \gamma$, $X(\sigma)$ is a simplicial complex of dimension
$l(\sigma)-1$. In particular, $X(\gamma)$ is a simplicial complex of dimension $n-1$.
\end{cor}

\begin{thm}\label{t:X = Z}
For $\theta_k \le \tau_i \le \tau_t$, we have
$c\left [X(\sigma , \rho) \right ] = Y(\sigma, \tau_i) = Z(\sigma, \tau_i)$.
\end{thm}

\emph{Proof:} It suffices to show that
$Z(\sigma, \tau_i) $ is contained in $c\left [X(\sigma , \tau_i) \right ]$ because
$c\left [X(\sigma , \rho) \right ]$ is contained in $Y(\sigma, \rho)$
(by definition) and
$Y(\sigma, \tau_i)$ is contained in $Z(\sigma, \tau_i)$
when $\tau_i \ge \theta_k$ (by theorem~\ref{t:tables}).
The proof is by induction on $i$, starting at the value $i_0$ for which
$\tau_{i_0} = \theta_k$.

\vskip .2cm

\textbf{Base step: }  Let $F_0= X(\sigma,\tau_{i_0-1})$
 be the
subcomplex of $X(\gamma)$ whose vertex set is
$\{\tau_1, \tau_2, \dots, \tau_{i_0-1}\}$.
Since $\langle \theta_1, \theta_2, \dots , \theta_{k} \rangle \in X(\sigma)$
(by corollary~\ref{c:firstsimplex}) and
$ \theta_1, \theta_2, \dots , \theta_{k-1} \in F_0$, it follows that
$\langle \theta_1, \theta_2, \dots , \theta_{k-1} \rangle \in F_0$.
Furthermore, if $\sigma_0 = R(\theta_k)\sigma$, then
$F_0 = X(\sigma_0, \tau_{i_0-1})$.   By induction on
$k = l(\sigma)$, we can assume that the assertion of the theorem
is valid if $\sigma$ is replaced by the length $k-1$ element $\sigma_0$ and hence
$c\left [X(\sigma_0, \tau_{i_0-1})\right ] =
Z(\sigma_0, \tau_{i_0-1})$.
(The base case of this inner induction is trivial since $k=1$ corresponds to a rank
$1$ group.)
Thus $|F_0|$ is convex and
$(k-2)$-dimensional.
\vskip .2cm
Now let $V_0= X(\sigma, \tau_{i_0})$ which has vertex set
$\{\tau_1, \tau_2, \dots, \tau_{i_0}\}$. Since $\mu(\theta_k)\cdot \tau_j=0$
whenever $\tau_j < \theta_k$ (as in the proof of corollary~\ref{c:firstsimplex}),
it follows that
$|V_0|$ is the cone with base $|F_0|$ and apex $\tau_{i_0}$.
Thus $|V_0|$ is convex and $(k-1)$-dimensional.
\vskip .2cm

The containment of $Z(\sigma,\tau_{i_0})$ in
$c\left [X(\sigma, \tau_{i_0})\right ]$
is demonstrated by examining the supports of the facets of the
positive cone, $c\left [ V_0 \right ]$, on $|V_0|$. Each support is of the form
$M(\sigma)\cap \mu(\tau_j)^{\perp}$ for some
$\tau_j$.

One of the facets  of $c\left [ V_0 \right ]$ contains $|F_0|$
and hence it has support $M(\sigma) \cap \mu(\tau_{i_0})^\perp$.
Each of the other facets
 of $c\left [ V_0 \right ]$ contains the vertex $\tau_{i_0}=\theta_k$ and
 hence its support is of the form
$M(\sigma)\cap \mu(\tau_j)^{\perp}$ for some $j \neq i_0$.
In fact, $\tau_j$ must belong to
$\{\theta_1 ,\dots  , \theta_{k-1}\} \cup \{\tau_j \mid \tau_j > \theta_k\}$
in this case.  For, $\mu(\tau_j)^\perp$
 cannot separate $|F_0|$ (because $|F_0| \subset |V_0| $ which is contained
 in one of the half-spaces bounded by $\mu(\tau_j)^\perp$)
 and hence if $\tau_j \notin \{\theta_1 ,\dots  , \theta_{k-1}\}$ then
 proposition~\ref{p:separate} implies that we must have $\tau_j > \theta_k$.
Thus the set of facets of
$c\left [ V_0 \right ]$ is of the form
\[ \{M(\sigma)\cap \mu(\tau_j)^{\perp} \mid
\tau_j \in \{\theta_{i_1}, \dots , \theta_{i_a}, \theta_k\} \cup
\{\tau_{j_1}, \dots, \tau_{j_b}\}\ \} \]
 where $j_l> i_0$ for $ l = 1,2, \dots , b$.
 It follows  that
 \begin{eqnarray*}
 Z(\sigma, \tau_{i_0}) &=& M(\sigma) \cap \mu(\theta_1)^{+}\cap \cdots
 \cap \mu(\theta_k)^{+}\cap \mu(\tau_{i_0+1})^-\cap \cdots \cap
 \mu(\tau_t)^-\\
 &\subseteq& M(\sigma) \cap \mu(\theta_{i_1})^{+}\cap \cdots
 \cap \mu(\theta_{i_a})^{+}\cap \mu(\theta_k)^+ \cap \\
 && \ \ \ \ \ \  \mu(\tau_{j_1})^-\cap \cdots \cap \mu(\tau_{j_b})^-\\
 &=& c\left [ X(\sigma, \tau_{i_0}) \right ].
 \end{eqnarray*}
\vskip .2cm
\textbf{Inductive step: }  Assume now that $i \ge i_0$ and that
$c\left [ X(\sigma, \tau_{i})\right ] = Z(\sigma, \tau_{i})$.
Let $F$  and $V$ be the subcomplexes of $X(\gamma)$
whose sets of vertices are
\[\{\tau_j \mid 1 \le j \le i \mbox{ and }
\mu(\tau_{i+1})\cdot \tau_j = 0 \}\]
and
\[\{\tau_{i+1}\}\cup \{\tau_j \mid 1 \le j \le i \mbox{ and }
\mu(\tau_{i+1})\cdot \tau_j = 0 \}\]
respectively. Then $|V|$ is a cone with base $|F|$ and apex $\tau_{i+1}$.
 We prove that the closure,~$Z$, of
$Z(\sigma,\tau_{i+1}) \setminus Z(\sigma,\tau_{i})$ is contained in
the positive cone, $c\left [ V \right ]$, on $|V|$.
Since $|V|$ is contained in $|X(\sigma, \tau_{i+1})|$, it will then follow
 that $Z(\sigma,\tau_{i+1}) = Z \cup Z(\sigma,\tau_i)$ is
contained in $c\left [ X(\sigma, \tau_{i+1}) \right ]$, as required.

\vskip 0.2cm

First we show that  $F$ is
$(k-2)$-dimensional.
Denote the length~$(k-1)$ element $ R(\tau_{i+1})\sigma$ by $\sigma'$.
Note that $\sigma' \le R(\tau_{i+1})\gamma$ so that
$M(\sigma') \subset \mu(\tau_{i+1})^\perp$.
Apply the procedures of sections~\ref{s:subcomplex} and \ref{s:walls}
to $\sigma'$  (i) to obtain a simple system
$\{\delta'_1, \delta'_2, \dots , \delta'_{k-1}\}$ for the
 set $P'$ of positive
roots in $M(\sigma')$ and (ii)  to calculate the reordering
$\{\theta'_1, \theta'_2, \dots , \theta'_{k-1}\}$ of the set
$\{\epsilon'_1, \epsilon'_2, \dots , \epsilon'_{k-1}\}$ for which
\[ \sigma' = R(\theta'_{k-1})\dots R(\theta'_2)R(\theta'_1)\] where
$\theta_1' < \theta_2' < \ldots < \theta_{k-1}'$ and where $\epsilon_j'$ is
given by
\[\epsilon'_j = R(\delta'_1)R(\delta'_2)\dots R(\delta'_{j-1})\delta'_j
\mbox{\ \ for } j = 1,\dots,k-1.\]

\vskip 0.2cm

Since $\{\theta'_1, \theta'_2, \dots , \theta'_{k-1}\}
\subset \mu(\tau_{i+1})^\perp$, if we show that
$\tau_{i+1} > \theta_j'$ for $j=1,\dots,k-1$, then it will follow that
$\langle \theta'_1, \theta'_2, \dots , \theta'_{k-1}\rangle \in F$
and hence $F$ is $(k-2)$-dimensional.

\vskip .2cm
Fix $j \in \{1,\dots,k-1\}$. From $\epsilon'_j \cdot  \mu(\tau_{i+1}) = 0$
we deduce that $\epsilon'_j \ne \tau_{i+1}$
and that
$R(\epsilon'_j) \le R(\tau_{i+1})\gamma$.  Therefore the length two element
$\sigma'' = R(\tau_{i+1})R(\epsilon'_j)$ precedes $\gamma$ and hence
$\sigma$ by equation~(\ref{e:induce}).

Assume now that $\epsilon'_j > \tau_{i+1}$ and let $Q = P \cap M(\sigma'')$.
Then, as in the proof of lemma~\ref{l:simpleaction},
 $\{\tau_{i+1},\epsilon_j'\}$ is the simple system for $Q$, with $\tau_{i+1}$
 the first root and $\epsilon_j'$ the last root. Thus
$R(\tau_{i+1})\epsilon'_j$ is a positive root and
$R(\tau_{i+1})\epsilon'_j > \tau_{i+1}$ (by lemma~\ref{l:simpleaction}).

Consider now $(\sigma')^{-1}(\epsilon_j')$ which is negative,
 by proposition~\ref{p:sigmaaction}. However
 \[(\sigma')^{-1} \epsilon'_j = (R(\tau_{i+1})\sigma)^{-1} \epsilon'_j
= \sigma^{-1}( R(\tau_{i+1})\epsilon'_j).\]
Since $R(\tau_{i+1})\epsilon'_j$ is positive, proposition~\ref{p:sigmaaction}
implies that $R(\tau_{i+1})\epsilon'_j= \epsilon_a$ for some $a$.
Thus $R(\tau_{i+1})\epsilon'_j$ precedes $\tau_{i+1}$ (since each $\epsilon_b$
does), which contradicts the earlier conclusion that
$R(\tau_{i+1})\epsilon'_j > \tau_{i+1}$. Thus the
assumption that $\epsilon_j'>\tau_{i+1}$ must have been false.

\vskip .2cm
Next we show that $|F|$ is convex. Since
$c\left [ X(\sigma, \tau_i) \right ] = Z(\sigma, \tau_i)$ by the inductive hypothesis,
we obtain
\begin{eqnarray*}
|F| &=& M(R(\tau_{i+1})\sigma) \cap c\left [ X(\sigma, \tau_i)\right ] \cap S^{n-1}\\
&=&  M(R(\tau_{i+1})\sigma) \cap Z(\sigma, \tau_i)\cap S^{n-1}
\end{eqnarray*}
which is convex.
\vskip .2cm
Now $|V|$, being a cone with a convex $(k-2)$-dimensional base,
must itself be convex and $(k-1)$-dimensional.
The proof that $Z$ is contained in $c\left [ V \right ]$ involves a close examination
of the facets of $c\left [ V \right ]$.  First one argues
that each facet of $c\left [ V \right ]$
is of the form~$M(\sigma) \cap \mu(\tau_j)^\perp$ for some
\[\tau_j \in \{\theta_{i_1}, \dots , \theta_{i_a}\} \cup
\{\tau_{i+1}, \tau_{j_1}, \dots, \tau_{j_b}\},\]
where $j_l > i+1$ for $l = 1,\ldots,b$. (This step is similar
to the corresponding argument for $c\left [ V_0 \right ]$,
with $\tau_{i+1}$ taking the
 place of $\theta_k$.) It then follows that
\begin{eqnarray*}
 Z &=& M(\sigma) \cap \mu(\theta_1)^{+}\cap \dots
 \cap \mu(\theta_k)^{+}\cap \mu(\tau_{i+1})^+\cap
 \mu(\tau_{i+2})^-\cap \dots \cap \mu(\tau_t)^-\\
 &\subseteq& M(\sigma) \cap \mu(\theta_{i_1})^{+}\cap \dots
 \cap \mu(\theta_{i_a})^{+}\cap \mu(\tau_{i+1})^+ \cap
 \mu(\tau_{j_i})^-\cap \dots \cap \mu(\tau_{j_b})^-\\
 &=& c\left [ V \right ] \ \ \ \ \ \mbox{as required.}
 \end{eqnarray*}
\vskip .2cm

\begin{cor}\label{c:char}
For each $\sigma \le \gamma$ the set
$|X(\sigma)|$ is convex.
Furthermore, $|X(\sigma)|$ is the intersection with $S^{n-1}$ of the
positive cone on the set $P$ of positive roots whose reflections
precede $\sigma$ and
\[|X(\sigma)| = S^{n-1} \cap M(\sigma) \cap \mu(\theta_1)^+ \cap \dots
\cap \mu(\theta_k)^+.\]
\end{cor}

\emph{Proof:} Apply theorem~\ref{t:X = Z} with $i = t$.

\vskip .2cm
We are now in a position to prove that $[I, \gamma]$ is a lattice.
\begin{thm}\label{t:lattice}
If $W$ is a finite real reflection group equipped
with the partial order $\le$ defined by reflection length and $\gamma$
is a Coxeter element for $W$, then the interval
$[I, \gamma]$ is a lattice.
\end{thm}

\emph{Proof:}  Choose a simple system and Coxeter element~$\gamma$
for $W$ as in section~\ref{s:root-mu}.
For $\sigma \leq \gamma$ and $\sigma \neq I$, construct the simplicial
complex $X(\sigma)$ and
define $X(I)$ to be the empty set~$\emptyset$.
We have seen that  $X(\sigma)$  has dimension $l(\sigma)-1$
and its vertex set is
\[P_{\sigma} = \{\rho_i \mid 1 \le i \le nh/2
\mbox{ and } R(\rho_i) \le \sigma  \}.\]
Furthermore, $X(\sigma)$ is a subcomplex of $X(\gamma)$ and
by corollary~\ref{c:char}, $|X(\sigma)|$
is convex.

Suppose now that $ \alpha$ and $\beta$ both precede $\gamma$ and consider
the sub-complex $X(\alpha) \cap X(\beta)$ of $X(\gamma)$.
Since  $|X(\alpha) \cap X(\beta)| = |X(\alpha)| \cap |X(\beta)|$,
the set $|X(\alpha) \cap X(\beta)|$ is
convex by corollary~\ref{c:char}.  If $d$ denotes the dimension of $X(\alpha) \cap X(\beta)$ then
 $|X(\alpha) \cap X(\beta)|$ is a union of
$d$-dimensional simplices and, for each $d$-simplex
 $\langle v_0, \dots , v_d \rangle$ in $ X(\alpha) \cap X(\beta)$, the linear subspaces
$\mathrm{span}(\{v_0, \dots , v_d\})$ and $\mathrm{span}(X)$ coincide.
We  associate  an element $ \sigma \le \gamma$ to $X(\alpha) \cap X(\beta)$
as follows.  Choose a $d$-simplex
$\langle v_0, \dots , v_d \rangle \in X$ and assume (without loss of generality)
that $ v_0 \le v_1 \le \dots \le v_d$.
Then
\[R(v_0)R(v_1)\dots R(v_d) \le \gamma^{-1}  \mbox{ \ by  definition of } X(\gamma)\]
and we define $\sigma$ to be $R(v_d)\dots R(v_1)R(v_0)$.  Since
$ \mathrm{span}(v_0,\dots,v_d)= M(\sigma) = \mathrm{span}(X(\alpha) \cap X(\beta))$,
it follows that $M(\sigma)$ contains and is spanned by the vertex set of
$X(\alpha) \cap X(\beta)$.  However this vertex set is $P_{\alpha} \cap P_{\beta}$ and
the theorem now follows by our remarks in section~\ref{s:poset}.

\section{Relationship with generalised associahedra.}\label{s:fomin}
In this section we embed the complex $X(\gamma)$ in a larger
simplicial complex $EX(\gamma)$ whose vertex set consists of all
positive roots and the negatives of the simple roots, and we show
that the geometric realisation of $EX(\gamma)$ is a sphere. If
$W$ is crystallographic, then we show that $EX(\gamma)$ is
simplicially isomorphic to the simplicial generalised associahedron for $W$
(defined in \cite{fz}).

Recall that we have partitioned the simple roots into two commuting sets,
$S_1=\{\alpha_1,\ldots,\alpha_s\}$ and
$S_2=\{\alpha_{s+1},\dots,\alpha_n\}$.  We will use the notation
\[-S_1 = \{-\alpha_1, \dots , -\alpha_s\} \,\,\, \mbox{ and }
\,\,\, -S_2 = \{-\alpha_{s+1}, \dots , -\alpha_{n}\}.\]
We note that the subscripting on
the $\rho$'s can be applied to negative indices with the convention that
$\rho_{-k} = \rho_{nh-k}$.
\begin{defn}\label{d:EX}
We define a set, $EX =EX(\gamma)$, of simplices
by declaring that
\begin{itemize}
\item
the vertex set is the ordered set
\[\{\rho_{-n+s+1}, \dots, \rho_{-1}, \rho_{0}\}
\cup \{\rho_1, \rho_2, \dots , \rho_{nh/2}\}
\cup \{\rho_{nh/2+1}, \dots , \rho_{nh/2+s}\},\]
\item
that an edge joins $\rho_i$ to $\rho_j$ whenever $i < j$,
$\rho_i \ne -\rho_j$  and $R(\rho_i)R(\rho_j)
\le \gamma^{-1}$ and
\item
$\langle \rho_{i_1}, \rho_{i_2}, \dots , \rho_{i_k}\rangle$
forms a $(k-1)$-simplex
if the vertices are distinct and pairwise joined by edges.
\end{itemize}
\end{defn}

We note that it follows from the definitions that the simplices
of $X(\gamma)$ are simplices of $EX(\gamma)$.  The extra vertices
are precisely the negatives of the simple roots $\{\alpha_1, \dots,
\alpha_n\}$.  Firstly,
 for $k = 1, \dots , n-s$  $\rho_{s+k} = -\gamma(\alpha_{s+k})$
implies that $\rho_{-n+s+k} = -\alpha_{s+k}$.  Secondly the set
$\{\rho_{nh/2+1}, \dots , \rho_{nh/2+s}\}$ is a permutation of
$\{-\alpha_1, \dots, -\alpha_s\}$.

\begin{thm}\label{t:crosspoly}
$EX(\gamma)$ is a simplicial complex and
$|EX(\gamma)|$ is a sphere of dimension $n-1$.
\end{thm}

\emph{Proof:}
Let $C$ be the spherical cross-polytope whose  set
of vertices is
$\{\pm \alpha_1, \dots , \pm \alpha_n\}.$
%, -\alpha_1, \dots , \dots , -\alpha_n\}\]
Then $C$ is a simplicial complex whose geometric realisation~$|C|$
is the  unit sphere in $\mathbf{R}^n$.
We will show that $EX(\gamma)$ is a simplicial
subdivision of $C$.

 \vskip .2cm
 Consider the  simplicial subdivision $C'$ of $C$ which is defined as
 follows.

 Let $K$ be the sub-complex of $C$ which consists
 of the simplex $A=\langle\alpha_1, \dots ,  \alpha_n\rangle$
  and all its faces.
We proved earlier that $|X(\gamma)| =|K|$. Thus  $X(\gamma)$ is a
simplicial subdivision of $K$.

\vskip 0.2cm

Extend this subdivision to the rest of $C$ as follows.
Any simplex of $C$ whose vertices are contained in
$\{- \alpha_1, \dots , - \alpha_n\}$  is not subdivided.
If $B$ is a simplex of $C$ which contains both positive and negative
roots then we can
write
\[ B= \langle\alpha_{i_1}, \dots,\alpha_{i_p} \rangle *
\langle -\alpha_{j_1}, \dots, -\alpha_{j_q} \rangle \]
where $i_1<\cdots<i_p$, $j_1<\cdots<j_q$, $*$ denotes the spherical join,
and $\{i_1,\dots,i_p\}
\cap\{j_1,\dots,j_p\} = \emptyset$. We extend the subdivision of
the simplex $\langle\alpha_{i_1}, \dots,\alpha_{i_p} \rangle$
(in $X(\gamma)$)
by taking the join of this subdivision with the simplex
$\langle -\alpha_{j_1}, \dots, -\alpha_{j_q} \rangle $.

\vskip 0.2cm

Suppose $J \subset \{1, \dots, s\}$, $K \subset \{s+1, \dots, n\}$
and that $\{\rho_{i_1}, \dots, \rho_{i_a}\} $ is an ordered set of
positive roots. Then $\{-\alpha_j \mid j \in J\cup K\}
\cup\{\rho_{i_1}, \dots, \rho_{i_a}\}$ is the vertex set for a
simplex in $C'$ or in $EX(\gamma)$ if and only if
\[R(\rho_{i_1})\dots R(\rho_{i_a}) \le
\left(\prod_{l\in K}R(\rho_l)\right)\gamma^{-1} \left(\prod_{m\in
J}R(\rho_m)\right).  \]
It now follows that $C' = EX(\gamma)$ as required.

\vskip 0.5cm

For convenience we recall some facts from~\cite{fz}
and express them in a manner consistent with our earlier notation.
If $W$ is a crystallographic finite reflection group
then the simplicial generalised associahedron, $GA(W)$, for the
simple system~$S_1 \cup S_2$
is a simplicial complex whose set of vertices, denoted $\Omega_{\ge -1}$
consists of all the positive roots and the negative simple roots.
Two piecewise-linear involutions $\tau_{+}$ and $\tau_{-}$ are
introduced in \cite{fz}. It can be shown that they are determined
by
\[\tau_+(\beta) =
\left \{\begin{array}{cc}
  R_1R_2\dots R_s(\beta) & \mbox{ if } \beta \not \in -S_2 \\
  \beta & \mbox{ if } \beta \in -S_2, \\
\end{array}
\right.\]
\[\tau_-(\beta) =
\left \{\begin{array}{cc}
  R_{s+1}R_{s+2}\dots R_n(\beta) & \mbox{ if } \beta \not \in -S_1 \\
  \beta & \mbox{ if } \beta \in -S_1. \\
\end{array}
\right.\] From this  we deduce that the action of $\tau_+\tau_-$
on $\Omega_{\ge -1}$ is given by
\[\tau_+\tau_-(\beta) =
\left \{\begin{array}{cc}
  \gamma (\beta) & \mbox{ if } \beta \not \in (-S_1) \cup S_2 \\
  -\beta & \mbox{ if } \beta \in (-S_1) \cup S_2 , \\
\end{array}
\right.\]
and the action of the inverse $\tau_-\tau_+$ is given by
\[\tau_-\tau_+(\beta) =
\left \{\begin{array}{cc}
  \gamma^{-1}(\beta) & \mbox{ if } \beta \not \in (-S_2) \cup S_1\\
  -\beta & \mbox{ if } \beta \in (-S_2) \cup S_1. \\
\end{array}
\right.\]

The compatibility degree  $(\alpha \parallel \beta )$ of any two
elements $\alpha,\beta \in \Omega_{\ge -1}$ is characterised in
\cite{fz} by the conditions
\begin{description}
\item[\it (i)] $(-\alpha_i \parallel \alpha ) = \max \{[\alpha :
\alpha_i], 0\}$, where $[\alpha:\alpha_i]=\alpha\cdot\mu_i$ for
$i=1,\dots,n$, and \item[\it (ii)] $(\alpha
\parallel \beta) = (\tau_{\pm}\alpha \parallel \tau_{\pm}\beta)$
for all $\alpha$ and $\beta$. \end{description}

We recall that two vertices $\alpha$ and $\beta$ in $GA(W)$ are
connected by an edge if and only if $(\alpha
\parallel \beta ) = 0$.
Finally, we note that $GA(W)$ is completely determined by its
one-skeleton.

\begin{thm}\label{t:genass}
If $W$ is crystallographic, then
$EX(\gamma)$ is simplicially isomorphic to
the  simplicial generalised associahedron $GA(W)$.
\end{thm}

\emph{Proof:} We continue to use the earlier notation.
First note that the vertex set $\Omega_{\ge -1}$ of
the associahedron is the same as that of $EX(\gamma)$, and we will refer
to the vertices in both complexes using the ordering from $EX(\gamma)$.
Since each of $EX(\gamma)$ and $GA(W)$ is
determined by its one-skeleton,  it suffices to show that
two vertices $\alpha$ and $\beta$ are joined by an edge in
$EX(\gamma)$ if and only if they are joined by an edge
in $GA(W)$.
To simplify notation we will write
$\rho_i \stackrel{GA}{\to} \rho_j$ or
$\rho_i \stackrel{EX}{\to} \rho_j$
if $\rho_i$ and $\rho_j$
are connected by an edge in the simplicial complex~$GA(W)$ or $EX(\gamma)$
respectively.
Note that the ordering on the vertices of $EX(\gamma)$ is such that they are arranged
into the following sets in the given order
\[-S_2, S_1, \gamma(-S_2),\gamma(S_1),\gamma^2(-S_2),
 \dots , \gamma^{-1}(-S_1) ,S_2, -S_1.\]

Suppose $\rho_i < \rho_j$.  Then $\rho_i \stackrel{EX}{\to} \rho_j$ is equivalent to
$R(\rho_i)R(\rho_j) \le \gamma^{-1}$.  First we can assume $\rho_j \ne -\rho_i$
since neither complex contains an edge from $\rho_i$ to $-\rho_i$.
\vskip .2cm
\textbf{Case 1.}   If $\rho_i \in -S_2$ then
$\rho_i = -\alpha_p$ for some $p$ satisfying $s+1 \le p \le n$.
Then $\rho_i \stackrel{GA}{\to} \rho_j$ if and only if $(-\alpha_p
\parallel \rho_j) = 0$.
However this is equivalent to $R(\rho_j)
\le R(\alpha_1)\dots R(\alpha_{p-1})R(\alpha_{p+1}) \dots
R(\alpha_n)$ by equation~(\ref{e:induce}).
But  $R(\alpha_1)\dots
R(\alpha_{p-1})R(\alpha_{p+1}) \dots R(\alpha_n)  = \gamma
R(\alpha_p)$ since $p > s$.
Thus $\rho_i \stackrel{GA}{\to}
\rho_j$ is equivalent by equation~(\ref{e:calc}) to
$R(\rho_j)R(\rho_i) \le \gamma$ which in
turn is the criterion for $\rho_i \stackrel{EX}{\to} \rho_j$.
\vskip .2cm
\textbf{Case 2.}   If $\rho_i \in S_1$, then $\rho_i =
\alpha_p$, for some $p$ satisfying $1 \le p \le s$ and
$\tau_-\tau_+(\rho_i) = -\rho_i = -\alpha_p \in -S_1$.
\vskip .2cm
\textbf{Case 2(a).}   It is possible that $\rho_j \in S_1$ also.
If this is the case then $\tau_-\tau_+(\rho_j) = -\rho_j \in -S_1$
also and applying $\tau_-\tau_+$ to both roots gives
\[(\rho_i \parallel \rho_j) = (-\rho_i \parallel -\rho_j) = 0\]
since $\rho_j$ and $\rho_i$ are distinct simple roots.
By definition of edges in the generalised associahedron we
 must have $\rho_i \stackrel{GA}{\to} \rho_j$.  However, in this subcase,
$\rho_j = \alpha_q$ with $1 \le p < q \le s$, so that
\[R(\rho_j)R(\rho_i) = R(\rho_i)R(\rho_j) \le \gamma \]
which is again equivalent to the criterion for
$\rho_i \stackrel{EX}{\to} \rho_j$.
\vskip .2cm
\textbf{Case 2(b).}  Now we assume $j > s$.  In this subcase,
we know that $\tau_-\tau_+(\rho_j) = \gamma^{-1} (\rho_j)$.
We have $\rho_i \stackrel{GA}{\to} \rho_j$ if and only
if $(-\alpha_p \parallel \gamma^{-1}(\rho_j)) = 0$.
However this is equivalent by equation~(\ref{e:induce}) to
\[R[\gamma^{-1}(\rho_j)] \le R(\alpha_1)\dots R(\alpha_{p-1})R(\alpha_{p+1})
\dots R(\alpha_n).\]
Since  $1 \le p \le s$,
we have $ R(\alpha_1)\dots R(\alpha_{p-1})R(\alpha_{p+1})
\dots R(\alpha_n)=  R(\alpha_p)\gamma$.  Thus $\rho_i \stackrel{GA}{\to} \rho_j$
is equivalent to
$R(\rho_i)R[\gamma^{-1}(\rho_j)] \le \gamma$ by equation~(\ref{e:calc}).  But
\[R(\rho_i)R[\gamma^{-1}(\rho_j)] \le \gamma
\Leftrightarrow R(\rho_i) \le \gamma R[\gamma^{-1}(\rho_j)] =
R(\rho_j)\gamma,\]
using equations~(\ref{e:calc}) and (\ref{e:conj}).  The last
condition in turn is the criterion for
$\rho_i \stackrel{EX}{\to} \rho_j$.
\vskip .2cm
\textbf{Case 3.}  If $i > s$ let $m$ be the smallest positive integer
with the property that
\[(\tau_-\tau_+)^m(\rho_i) \in (-S_2)\cup S_1.\]
Now $\rho_i \stackrel{GA}{\to} \rho_j$ if and only if
$(\tau_-\tau_+)^m(\rho_i)\stackrel{GA}{\to} (\tau_-\tau_+)^m(\rho_j)$.
However, Case 1 or Case 2 now applies to this new pair of roots and
$(\tau_-\tau_+)^m = \gamma^{-m}$ when applied to $\rho_i$ and $\rho_j$.
Thus $\rho_i \stackrel{GA}{\to} \rho_j$ if and only if
\[R(\gamma^{-m}(\rho_j))R(\gamma^{-m}(\rho_i)) \le \gamma,\]
by the proof in cases 1 and 2.
But this last condition is equivalent to
$R(\gamma^{-m}(\rho_i)) \le R(\gamma^{-m}(\rho_j))\gamma$, or
$\rho_{i-mn} \cdot \mu_{j-mn} = 0$.
Since $\gamma$ is an isometry this is equivalent to $\rho_{i} \cdot \mu_{j} = 0$
and hence to $\rho_i \stackrel{EX}{\to} \rho_j$.

\begin{note}\label{n:noncrys}
We observe that theorem~\ref{t:genass} provides a new proof that the simplicial
generalised associahedron is a simplicial complex whose geometric
realisation is a sphere.  This proof is independent of the classification
of finite real reflection groups and extends the work of Fomin and Zelevinsky~\cite{fz}
 to include the
non-crystallographic finite reflection groups.
\end{note}

\end{document}